\documentclass[12pt]{amsart}
\usepackage{amsfonts}
\usepackage{amscd}
\usepackage{amsmath, amssymb}
\usepackage[all]{xy}
\CompileMatrices
\theoremstyle{plain} \numberwithin{equation}{section}
\newtheorem{theorem}{Theorem}[section]
\newtheorem{corollary}[theorem]{Corollary}

\newtheorem{lemma}[theorem]{Lemma}
\newtheorem{proposition}[theorem]{Proposition}
\theoremstyle{definition}

\newtheorem{remark}[theorem]{Remark}
\newtheorem{example}[theorem]{Example}

\newcommand{\one}{\mathbf{1}}
\newcommand{\COSH} {{\bf cosh}}
\newcommand{\COS} {\mathbf{cos}}
\newcommand{\la}{\langle}
\newcommand{\ra}{\rangle}
\newcommand{\cA}{\mathcal {A}}
\newcommand{\cC}{\mathcal{C}}
\newcommand{\cD}{\mathcal{D}}
\newcommand{\cF}{\mathcal{F}}
\newcommand{\cO}{\mathcal{O}}
\newcommand{\cH}{\mathcal{H}}
\newcommand{\cK}{\mathcal{K}}
\newcommand{\cR}{\mathcal{R}}
\newcommand{\cT}{\mathcal{T}}
\newcommand{\cW}{\mathcal{W}}

\newcommand{\Sl}{\mathrm{Sl}}
\newcommand{\SO}{\mathrm{SO}}
\newcommand{\ad}{\mathrm{ad}}
\newcommand{\Ad}{\mathrm{Ad}}

\newcommand{\bs}{\backslash}

\newcommand{\fa}{\mathfrak{a}}

\newcommand{\bc}{\mathbf{c}}

\newcommand{\fg}{\mathfrak{g}}
\newcommand{\fh}{\mathfrak{h}}
\newcommand{\fk}{\mathfrak{k}}

\newcommand{\fp}{\mathfrak{p}}
\newcommand{\fq}{\mathfrak{q}}

\newcommand{\ga}{\alpha}

\newcommand{\gO}{\Omega}
\newcommand{\cL}{\mathcal{L}}
\newcommand{\cS}{\mathcal{S}}

\newcommand{\D}{\mathbb{D}}

\newcommand{\R}{\mathbb{R}}
\newcommand{\C}{\mathbb{C}}

\newcommand{\N}{\mathbb{N}}
\newcommand{\Z}{\mathbb{Z}}
\newcommand{\bO}{{\bf{O}}}

\newcommand{\x}{{\bf x}}
\newcommand{\y}{{\bf y}}
\newcommand{\z}{{\bf z}}

\newcommand{\ZH}{Z_H}

\newcommand{\gOH}{\Omega_H}

\newcommand{\zH}{z_H}

\newcommand{\Ups}{\mathbf{\mathrm{Hor}}}

\newcommand{\Int}{\mathop{\textrm{int}}}
\newcommand{\im}{\mathop{\textrm{Im}}}

\newcommand{\Tr}{\mathop{\mathrm{Tr}} }

\newcommand{\eps}{\varepsilon}

\newcommand{\oline}{\overline}

\def\sideremark#1{\ifvmode\leavevmode\fi\vadjust{\vbox to0pt{\vss
 \hbox to 0pt{\hskip\hsize\hskip1em
 \vbox{\hsize2cm\tiny\raggedright\pretolerance10000
 \noindent #1\hfill}\hss}\vbox to8pt{\vfil}\vss}}}%

\theoremstyle{plain}

\begin{document}
\title[Holomorphic horospherical transform]{Holomorphic
horospherical transform on non-compactly causal 
spaces}

\author{Simon Gindikin, Bernhard Kr\"otz and Gestur \'{O}lafsson}
\address{Department of Mathematics, Rutgers University, New Brunswick, NJ 08903}
\email{gindikin@math.rutgers.edu}
\address{Max-Planck-Institut f\"ur Mathematik, Vivatsgasse 7, D-53111 Bonn,
Germany}
\email{kroetz@mpim-bonn.mpg.de}
\address{Department of Mathematics, Louisiana State
University, Baton Rouge,
LA 70803, USA}
\email{olafsson@math.lsu.edu}


\thanks{SG was supported in part by NSF-grant DMS-0070816}
\thanks{BK was supported by the RiP-program in Oberwolfach,
and a Heisenberg fellowship of the DFG}
\thanks{G\'O was supported by the RiP-program in Oberwolfach,
NSF grants DMS-0139783 and DMS-0402068}

\begin{abstract}We develop integral geometry for non-compactly 
causal symmetric spaces. 
We define a complex horospherical 
transform and, for some cases, identify it with a Cauchy type integral.
\end{abstract}
\maketitle

\section*{Introduction}

Within the class of semisimple symmetric spaces $Y=G/H$ we focus on those 
which can be realized as Shilov boundaries 
of Stein tubes $D=D(Y)$ in the affine complexification 
$Y_\C=G_\C/H_\C$. These $Y$ come in two different flavors:
compactly causal (CC) and non-compactly causal (NCC) symmetric spaces. 
It is important to mention that one can realize 
different series of representations of $L^2(Y)$ 
as holomorphic functions on $D$: 
the holomorphic discrete series for CC
and a multiplicity one subspace of the most continuous spectrum in the
NCC-case.

\par
In \cite{GKO} we developed integral geometry for $D$ in the CC-case.
If we consider the usual (real) horospherical transform on $Y$,
then holomorphic discrete series lie in its kernel.
So we considered a complex version of such a transform - horospherical Cauchy transform - using  a kernel of Cauchy type with singularities on
complex horospheres (on $Y_\C$) which do not intersect $Y$.
As a result we constructed a dual domain $\Xi_+$ in the manifold $\Xi$ of complex horospheres on $Y_\C$ and
our horospherical transform is an intertwining operator from holomorphic functions on
$D$ to holomorphic functions on $\Xi_+$ which admits an explicit inversion.

\par In this paper we develop holomorphic integral geometry for
NCC-spaces. In \cite{GKO} we constructed for any NCC-space a $G$-invariant tube
type domain $D(Y)$ in $Y_\C$ which has $Y$ as Shilov boundary. 
The crucial point is that $D$ contains the associated Riemannian
symmetric space $X=G/K$ as a totally real submanifold. 
We have $X_\C\simeq Y_\C$. Let us remark 
that all spherical functions on $X$
extend holomorphically to $D$ and admit 
boundary values on $Y$. 
Furthermore, in some cases,  $D$ coincides with the 
maximal $G$-domain of holomorphy of spherical
functions on $X$ -- the complex crown of $X$.

\par There are substantial differences in our constructions of holomorphic
horospherical transform for CC and NCC-spaces. In the CC-case we 
used a Cauchy type integral operator on $Y$ with singularities on complex
horospheres which do not intersect $Y$ (horospherical Cauchy transform).
For NCC-spaces we are able to  construct similar operators in some  cases, 
including the important subclass of Cayley type spaces. In general, 
we give another construction. We use the fact that complex
horospheres which do not intersect $Y$
allow real forms in $D$ and we can
construct a complex horospherical transformation extending 
the real one. Let us  mention 
that this construction does not work in CC-case. 
If we have both type of
constructions the second one can be regarded 
as a  residue of first one.

\par  Let us decsribe our results in more detail. 
For a fixed choice of horospherical coordinates 
$X\simeq N\times A$ one parameterizes 
the set $\Xi_\R$ of real horospheres on $X$ by 
$G/MN$ (with $M=Z_K(A)$ as usual). For 
Schwartz-functions $f$ on $X$ one defines 
the real Radon transform 
$$\cR_\R : \cS(X)\to C^\infty(\Xi_\R), \ \ \cR(f)(\xi)=\int_\xi f\,.$$     
Our first goal is 
to give a holomorphic 
version of $\cR_\R$. 
For that we construct a certain  $G$-invariant submanifold 
$\Xi_+$ of the complex horospheres $\Xi=G_\C/ M_\C N_\C$ on $X_\C
=Y_\C$. The complex horospheres of $\Xi_+$ have the property 
that they do not intersect $Y$. Let us point out, 
that in  contrast to the CC-case, $\Xi_+$ is not longer open in $\Xi$;  
it is a $CR$-submanifold of complex dimension  equal to the rank of $Y$ (i.e. 
$\dim A$). 
Next we recall the most-continuous Hardy-space $\cH^2(D)$ on 
$D$ from \cite{GKO}. This is a Hilbert space 
of holomorphic functions on $D$ with an equivariant 
isometric boundary map $b: \cH^2(D)\to L_{\rm mc}^2(Y)$ into a
multiplicity-one subspace of the most-continuous spectrum of 
$Y$. We show that there is a 
natural dense $G$-subspace $\cH^2(D)_0$ of $\cH^2(D)$ whose elements
restrict to Schwartz-functions on $X$ and for which 
$\cR_\R$ allows a holomorphic $G$-equivariant extension 

$$\cR: \cH^2(D)_0\to CR^\infty(\Xi_+)\, .$$  
We call $\cR$ the holomorphic horospherical Radon transform.

\par The main result of this paper is an interpretation  
of a (slightly twisted) Radon transform $\cR_\chi$ 
as a Cauchy-type integral. For the highly symmetric subclass of 
Cayley-type spaces CT=CC$\cap$NCC we show in Theorem 5.3 that 
there is Cauchy-type kernel 
$\cK_\chi(,)\in CR_\chi^\infty(\Xi_+)\hat \otimes \cS(Y)$ such that 
\begin{equation} \label{II1} 
\cR_\chi(f)(\cdot)={1\over (2\pi)^n} \int_Y b(f)(y)\ \cK_\chi (\cdot , y) \ dy \in CR_\chi^\infty(\Xi_+)\, .
\end{equation}
Note that (\ref{II1}) is in analogy to our definition of the 
holomorphic horospherical Radon transform for the CC-case \cite{gko2}.
For the general case of an NCC-space we believe that
there is an integral representation of $\cR_{\chi}$ as in (\ref{II1})
but we are not aware of such. 
We consider the definition of a holomorphic horospherical Radon transform
in terms of  a Cauchy type integral to be most 
robust for further generalization to other symmetric spaces or different 
series of representations.

\par For the CC-case, following \cite{G3}, we found an inversion 
formula for the holomorphic horospherical 
transform of compelling beauty: 
$$f = (\cL \cR f)^\vee\, .$$
Here $\cL$ is a differential operator and 
$\phi\mapsto \phi^\vee$ is a dual transform \cite{gko2}
with regard to a double fibration which brings $D$ and $\Xi_+$ 
into duality. 
Let us remark that  the situation is different for the 
NCC-spaces -- for the basic example $Y=\Sl(2,\R)/ \SO(1,1)$ we show that 
$\cL$ is a pseudo differential  but not a differential operator 
(cf. Section 6).

\par
The paper is concluded with a geometric definition of the most continuous
Hardy space $\cH^2(D)$ which is of independent interest. 
\medskip 
\par {\it Achknowledgement}: We thank the referee for his patient pressure 
on more precision. It is due to him that the paper is 
now more readable.

\section{Horospheres on symmetric spaces of triangular type } \label{s-one}
\noindent
In \cite{GKO} we associated to
every NCC-symmetric space 
$Y=G/H$ a $G$-Stein manifold $D$   with the following properties:
\begin{enumerate}
\item The complex manifold $D$ has a natural $G$-realization in the
complexification
$Y_\C$ of $Y$;
\item The symmetric space $Y=G/H$ is $G$-isomorphic to the
distinguished (Shilov) boundary
of $D$.
\end{enumerate}
The objective of this section is to study the space
$\Xi=G_\C/M_\C N_\C$
of horospheres  in $Y_\C$ in relation to $D$.
In particular we will introduce a natural $G$-invariant
$CR$-manifold $\Xi_+\subset
\Xi$ whose elements have the properties that they do non intersect
the real space $Y$, i.e. have no \textit{real points}.

\subsection{Notation}\label{ss-one}

We begin with the construction of NCC-symmetric spaces from the 
complex geometric point of view.  
\par Let us denote by $\fg$ a simple non-compact real Lie algebra.
We fix a Cartan involution $\theta:\fg\to \fg$, write 
$\fg=\fk+\fp$ for the associated eigenspace decomposition and 
select a maximal abelian subspace $\fa\subset\fp$. 
Write $\Pi=\{\alpha_1,\ldots,\alpha_n\}$  for a basis of the restricted root system 
$\Sigma=\Sigma(\fg,\fa)$ and expand the highest root 
$\beta$ in terms of $\Pi$
$$\beta=k_1 \alpha_1 +\ldots + k_n \alpha_n \qquad (k_i\in \Z_{>0})\, .$$
One calls $\alpha_i$ {\it miniscule} if  $k_i=1$. 
Let us recall that miniscule elements exist if and only if 
$\Sigma$ is not of the type $E_8$, $F_4$ or $G_2$. 
For the remainder we fix a miniscule root $\alpha$ and define $Z\in \fa$ 
by the requirement $\alpha_i(Z)=\delta_{\alpha_i \alpha}$. 
Observe that $Z$ gives rise to a triangular decomposition
\begin{equation}\label{3g} \fg=\fg_{-1} + \fg_0 + \fg_1 \end{equation}
with $\fg_j=\{ X\in \fg\mid [Z,X]=j X\}$.  
Further the $3$-grading (\ref{3g}) defines an involution 
$\sigma$ of $\fg$ by $\sigma(X)=(-1)^j X $ for $X\in \fg_j$ and 
eigenspace decomposition 
$$\fg=\fg_0 + ( \fg_{-1}+\fg_{1})\, .$$
By construction $\sigma$ and $\theta$ commute and hence 
$\tau=\sigma\circ \theta$ defines an involution whose 
eigenspace decomposition shall be denoted  
$\fg=\fh+\fq$. The so-obtained symmetric pairs 
$(\fg, \fh)$ we call non-compactly causal (NCC). 
It is an elementary exercise to classify 
all  NCC-symmetric pairs. 
\par We denote by $\cW$ the Weyl group of $\Sigma$, 
set $Z_H={\pi\over 2} Z$ and define 
$\Omega_H\subseteq \fa$ by
\begin{equation}\label{de-OH}
\gOH=\Int  \{\hbox{convex hull of }\ \cW(\ZH)\}\, .
\end{equation}
Observe that  $\oline \Omega_H$ is a compact convex subset of
$\fa$ with extreme points $\cW (Z_H)$.

\medskip On the global level we fix a linear connected Lie group $G$ 
with Lie algebra $\fg$. We denote by 
$G_\C$ the corresponding linear connected complexification of $G$. 
Further 
we request 
that $\tau$ as well is complex linear extension 
$\tau_\C$ exponentiate to involutions on $G$, resp. $G_\C$, and we denote 
by $H$, resp. $H_\C$ their corresponding 
fixed point groups. In this way we obtain a totally real embedding 
$$Y\hookrightarrow Y_\C= G_\C/H_\C$$
of $Y$ in the Stein symmetric space $Y_\C$.
We refer to $Y$ as a non-compactly causal 
symmetric space (NCC).

\par Attached to $Y$ and $\Omega_H$ comes a
Stein manifold $D$ which we will now describe. 
Write
$K<G$ for the compact group of $\theta$-fixed elements and
$X=G/K$ for the corresponding Riemann symmetric space. As before
we obtain a totally real embedding

$$X\hookrightarrow X_\C= G_\C/K_\C\, .$$

Recall that $H_\C$ and $K_\C$ are conjugate, i.e. with  $\zH=\exp
(i\ZH)$ we have cf. \cite{GKO}:
\begin{equation}\label{eq=conj}
e^{i\ad (\ZH)}\fk_\C=\fh_\C\qquad \mathrm{and}\qquad \Ad (\zH)K_\C= H_\C\, .
\end{equation}
Hence $X_\C$ and $Y_\C$  are canonically $G_\C$-isomorphic
via the map
$$X_\C\ni gK_\C \mapsto  g\zH^{-1}  H_\C\in Y_\C\, .
$$
In the sequel we identify $X_\C$ with $Y_\C$.

\par We write $x_o=K_\C\in X_\C$ for the base point in  $X_\C$ and set
$$D=G\exp(i\Omega_H)\cdot x_o\, .$$
Note that $D$ was denoted by $\Xi_H$ in our previous article
\cite{GKO}.
According to \cite{GKO},  $D$
is an open $G$-invariant Stein neighborhood of $X$ in $X_\C=Y_\C$.
Moreover,  the map
$Y=G/H\ni gH\mapsto   gz_H\cdot x_o\in X_\C $ identifies $Y$ with
the distinguished
boundary $ \partial_dD$
of $D$ (see \cite{GKO}, Section 1, for more
details).
\par In summary, the symmetric Stein manifold $X_\C=Y_\C$
admits 2 real forms $X$ and $Y$ and a Stein neighborhood $D$ of $X$
with $Y$ as its Shilov boundary.

\subsection{Complex horospheres}
\label{ss-two} 
From now we let $Y=G/H$ be an NCC symmetric space and $X=G/K$ its 
Riemannian counterpart. 
\par In this section we introduce the $G$-space of
\textit{horospheres} in the complex manifold $X_\C=Y_\C$. This was done
for CC-spaces in \cite{gko2}.

\par We begin with some general remarks on convexity which we will
use frequently. Let $G=NAK$ be an Iwasawa decomposition of $G$ and
$N_\C A_\C K_\C\subsetneq G_\C$ its Zariski-open complexification.
In particular, $N_\C A_\C\cdot x_o$ is a Zariski-open subset in the
affine variety $X_\C$. Define the finite $2$-group $F=A_\C \cap
K_\C$ and  note that there are well defined holomorphic maps

$$n: N_\C A_\C \cdot x_o\to N_\C, \qquad a: N_\C A_\C\cdot x_o \to A_\C/F$$
such that $z=n(z) a(z)\cdot x_o$ for all $z\in N_\C A_\C\cdot x_o$.
Now, the fact that $D$ is contractible and $D\subset N_\C A_\C\cdot
x_o$, implies that $a|_{D}$ admits a well defined holomorphic
logarithm
$$\log a: D \to \fa_\C\, .$$
For $Z\in \Omega_H$, the complex convexity theorem (cf.\
\cite{GK02,ko}) then implies that
\begin{equation} \label{cc} \im \log  a(G\exp (iZ)\cdot x_o)={\rm conv} (\cW\cdot Z)
\end{equation}
where ${\rm conv} (\cdot)$ denotes the convex hull of $(\cdot)$.

\par Submanifolds of $X_\C$ of the type
$$gN_\C \cdot x_o \qquad (g\in G_\C)$$
will be referred as {\it horospheres}. We denote by $\Ups(X_\C)$ the set of all
horospheres on $X_\C$ and note that
$\Ups (X_\C)$ has a natural $G$-structure
$(g,hN_\C\cdot x_o)\mapsto gh N_\C\cdot x_o$.

\par To understand the space horospheres and
the related harmonic analysis it is useful to bring them in the
context of  a double fibration. Set $M=Z_K(A)\subset M_\C=Z_{K_\C}
(A)$, define
$$\Xi=G_\C/M_\C N_\C$$
and consider:

\begin{equation}\label{eq=df}
 \xymatrix { & G_\C/ M_\C \ar[dl]_{\pi_1} \ar[dr]^{\pi_2} &\\
\Xi & & X_\C\,.}\end{equation}
Then horospheres in $X_\C$
are exactly the subsets of $X_\C$ of the
form
\begin{equation}\label{eq=E} E(\xi)=\pi_2(\pi_1^{-1}(\xi)) \qquad (\xi\in \Xi)\, .
\end{equation}
If $\xi_o=M_\C N_\C\in \Xi$ denotes the base point and $\xi=g\cdot
\xi_o\in \Xi$ then, using that $M_\C \subset H_\C$, we have:
$$E(\xi)=gM_\C N_\C\cdot x_o=gN_\C \cdot x_o\subset X_\C\, . $$

\par
Similarly,  for $z\in X_\C$ we set
\begin{equation}\label{eq=S} S(z)=\pi_1(\pi_2^{-1}(z))\ .\end{equation}
If $z=g\cdot x_o$ for $g\in G_\C$, then
$S(z)=gK_\C\cdot \xi_o$.
Moreover, for $z\in X_\C$ and $\xi\in \Xi$ one has
the incidence relations
\begin{equation}\label{eq=indi}
z\in E(\xi)\iff \pi_1^{-1}(\xi)\cap \pi_2^{-1}(z)\neq \emptyset
\iff \xi\in S(z)\, .\end{equation}

\begin{proposition}\label{gen-parametrization} The map
$$\Xi\to \Ups(X_\C), \ \
\xi\mapsto E(\xi)$$
is a $G_\C$-equivariant bijection.
\end{proposition}

\begin{proof}  $G_\C$-equivariance and
surjectivity are clear.
The injectivity follows the same way as in the proof of
Proposition 2.1 in \cite{gko2} by replacing $H_\C$ by
$K_\C$.\end{proof}

One of the important features of $\Xi$ is, that there exists a right
$A_\C$-action on $\Xi$ that commutes with the left $G_\C$-action.
For $\xi=g\cdot\xi_o$ and $a\in A_\C$ we set

\begin{equation}\label{eq=ra} \xi\cdot a= ga\cdot \xi_o\, .
\end{equation}
Since $A_\C$ normalizes $M_\C N_\C$ it is clear that (\ref{eq=ra})
is well defined. {}From the definition it is also clear that the
left $G_\C$-action and the right $A_\C$-action commutes. In this way
we obtain an action of $G_\C\times A_\C$ on $\Xi$ by

$$(G_\C\times A_\C) \times \Xi \to \Xi, \ \ ((g,a), \xi)\mapsto
g\cdot \xi\cdot a\, .$$

\par
We conclude this subsection with an alternative characterization of
horospheres as level sets of holomorphic functions.
For that let $\{ \omega_1, \ldots, \omega_n\}\subset \fa^*$ be
the set of fundamental $K$-spherical lowest weights. For each $1\leq j\leq n$ we write
$(\pi_j, V_j)$ for the corresponding finite dimensional representation of $G$
with lowest weight $\omega_j$. We extend this representation to a
holomorphic representation of $G_\C$ which we denote by the same symbol.
Endow $V_j$ with a complex bilinear pairing $ \la  \hbox to 0.5em{},  \hbox to 0.5em{} \ra$ such that
$\la \pi_j(g)v, w\ra =\la v, \pi_j(\theta(g)^{-1})w\ra$ for all
$v,w\in V_j$ and $g\in G_\C$. Such a form  $\la \hbox to 0.5em{} ,\hbox to 0.5em{} \ra$
exists as $\pi\circ \theta$ is isomorphic to the representation
contragredient to $\pi_j$.
We write $v_j\in V_j$ for a lowest weight
vector and $\eta_j\in V_j$ for a $K_\C$-fixed vector subject to the
normalization $\la \eta_j, v_j\ra=1$. Finally, define holomorphic functions
on  $f_j: G_\C \to \C$ by

\begin{equation} \label{ef} f_j(g)=\la \pi_j(g)\eta_j, v_j\ra \qquad (g\in G_\C)\,.
\end{equation}
Note, that we have
\begin{equation}\label{eq-NAKcov}
f_j(nak)=a^{\omega_j}
\end{equation}
for all $n\in N_\C$, $k\in K_\C$ and $a\in A_\C$. Here, as elsewhere in this
article, we use the notation
$a^\mu =e^{\mu (X)}$ if $a=\exp X\in A_\C$.
We recall that
\begin{equation}\label{eq-zero}
G_\C\setminus N_\C A_\C K_\C=\{g\in G_\C\mid \prod_{j=1}^n f_j(g)=0\}\, .
\end{equation}
(see \cite{vdB}, Lemma 3.4). 

\begin{lemma}\label{le-NK}

$N_\C K_\C=\{g\in G_\C \mid f_j (g)=1 \quad \hbox{for all}\quad 1\leq j\leq n\}$.
\end{lemma}
\begin{proof} This follows from (\ref{eq-NAKcov}) and
(\ref{eq-zero}).\end{proof}

We will often view $f_j$,
or more generally it left translates, as a
function on $X_\C$. We will also, without further comments, view
the function $g\mapsto f_j(g^{-1}z)$, $z\in X_\C$ as a function on
$\Xi$. With that in mind we have:

\begin{lemma} Let $\xi \in \Xi$ and
$x\in X_\C$. Then
$$E(\xi)=\{ z\in X_\C \mid f_j(\xi^{-1} z )=1 \quad \hbox{for all}
\quad 1\leq j\leq n\}$$
and
$$S(x)=\{ \varrho \in \Xi \mid f_j(\varrho^{-1} x )=1 \quad \hbox{for all}
\quad 1\leq j\leq n\}$$
\end{lemma}

\begin{proof} Notice that $E(g\cdot \xi_o)=gE(\xi_o)$ and
$S(g\cdot x_o)=gS(x_o)$. We can therefore assume that
$\xi =\xi_o$ and
$x=x_o$. Now, the claim is  a reformulation of
Lemma \ref{le-NK}.
\end{proof}

\subsection{Some $G$-submanifolds of $\Xi$}
We define
the $G$-space  of {\it real horospheres} in
$X$ as
$$\Xi_\R = G/MN\, .$$
Then $\Xi_\R\subset \Xi=G_\C / M_\C N_\C$ is
obviously a totally real $G$-submanifold of $\Xi$ and the right
$A$-action leaves $\Xi_\R$ invariant.
\par Let  $T=\exp(i\fa)\subset G_\C$ and note that $A_\C=A\times T$; note that
$F=K\cap T$. We contrast $\Xi_\R$ with the  $G\times A_\C$-invariant subset
of $\Xi$

\begin{equation} \Xi_0=G\cdot \xi_o\cdot A_\C =GA_\C\cdot \xi_o\, .
\end{equation}

\begin{proposition} The following assertion hold:
\begin{enumerate}
\item The map
\begin{equation} \label {eq=decomp} K/M\times_F  A_\C \to \Xi_0, \ \ [kM,a]\mapsto ka\cdot \xi_o
\end{equation}
is a real analytic isomorphism.
\item The map
$$\Xi_\R \times_F  T\to \Xi_0, \ \ [gMN, t]\mapsto gt\cdot\xi_o$$
is a $G$-equivariant real analytic diffeomorphism.
\end{enumerate}
\end{proposition}
\begin{proof} (i) follows from the fact that
$G=KAN$ and $N A_\C \subset A_\C N_\C$.
Finally (ii), is a consequence of (i).
\end{proof}

Note that (\ref{eq=decomp}) describes a natural $CR$-structure on $\Xi_0$
of $CR$-dimension $\dim A$ and $CR$-codimension $\dim K/M$.

Define a tube domain in $A_\C$ by
$$T (\Omega_H)=\exp (\fa + i\Omega_H)
=A\exp (i\Omega_H)\simeq \fa +i\Omega_H\, $$
and set
\begin{equation}\label{defYH}
\Xi_+=G\exp(i\Omega_H)\cdot \xi_o=K T(\Omega_H)\cdot \xi_o\, .
\end{equation}
Then $\Xi_+$ is a real-analytic, $G$-invariant open submanifold of $\Xi_0$.
In particular $\Xi_+$ is a $CR$-manifold.
The coordinate decomposition of $\Xi_0$ slightly simplifies
for $\Xi_+$.
\begin{proposition} For $\Xi_+$ the following assertions hold:
\begin{enumerate}
\item The map
$$K/M\times T(\Omega_H) \to \Xi_+, \ \ (kM,a)\mapsto ka\cdot \xi_o$$
is a real analytic isomorphism.
\item The map
$$\Xi_\R  \times \Omega_H \to \Xi_+, \ \ (gMN, Z)\mapsto g\exp(iZ)\cdot\xi_o$$
is a $G$-equivariant real analytic diffeomorphism.
\end{enumerate}
\end{proposition}

We conclude this section with a remark on
the structure of $\Xi_+$.

\begin{remark} (Shilov boundary of $\Xi_+$) The map
$$\Xi_\R\to \partial\Xi_+, \ \ gMN\mapsto gz_H\cdot \xi_o$$
identifies $\Xi_\R$ as the Shilov boundary $\partial_S \Xi_+$ of
$\Xi_+$. In this sense $\Xi_\R$ parameterizes the real horospheres on
$Y$ (see also Lemma \ref{nc} below).
\end{remark}

\subsection{Horospheres without real points}\label{ss=three}

The aim of this subsection is to show
that horospheres corresponding to $\Xi_+$ do not
contain real points, i.e., are disjoint form
$Y$.

Recall from Subsection \ref{ss-one} that we identify $Y=G/H$ with
the (Shilov) boundary orbit $G\cdot y_o\subset X_\C$ of $y_o=z_H\cdot x_o$ in
$G_\C/K_\C$. Define the parameter set of \textit{horospheres without
real points} by
\begin{equation}\label{de-Ur}
\Xi_{\rm nr}=\{ \xi \in \Xi\mid E(\xi) \cap Y=\emptyset\}\, .
\end{equation}
The following statement should be compared to  the complex
convexity Theorem (\ref{cc});it means that convexity breaks down at
the extreme points of $\Omega_H$.

\begin{lemma}\label{nc}
Let ${\mathcal U}=\bigcup_{w\in \cW} NAwH$. Then ${\mathcal U}\cdot y_o$ is
open and dense in $G\cdot y_o$ and
$$G\cdot y_o\cap N_\C A_\C \cdot x_o
= {\mathcal U}\cdot y_o=
N A \cW z_H\cdot x_o\, .$$
\end{lemma}

\begin{proof} It is a special case of a Theorem by
Rossmann-Matsuki (cf. \cite{Ma79}) that ${\mathcal U}$ is dense in
$G$. Hence ${\mathcal U}\cdot y_o$ is dense in $G\cdot y_o$. As
$\zH^{-1}H_\C \zH =K_\C$ (cf. \ (\ref{eq=conj})),  it follows that
${\mathcal U}\cdot y_o = NA\cW \zH\cdot x_o$. It remains to show the
inclusion "$\supseteq$" for the first asserted equality. But this
follows from (\ref{eq-zero}).
\end{proof}

We can now prove the main result of this subsection.
\begin{theorem} \label{ti}$\Xi_+\subseteq \Xi_{\rm nr}$.
\end{theorem}

\begin{proof} We argue by contradiction. Note that $E(\xi)\cap Y\neq \emptyset$ for some $\xi\in\Xi_+$ means that
there exist $Z\in \Omega_H$ such that
$$Gz_H\cdot x_o \cap \exp(iY)N_\C\cdot x_o\not= \emptyset\, .$$
Now  $\exp(iY)N_\C\cdot x_o=N_\C \exp(iY)\cdot x_o\subset N_\C \exp(i\Omega_H)\cdot x_o$
and the assertion follows from Lemma \ref{nc}.
\end{proof}

\subsection{Real forms of $E(\xi)$ and $S(z)$}
In this last part of this section we introduce certain $G$-invariant
real forms of the
complex manifolds $E(\xi)$ and $S(z)$.

\par We begin with the horospheres. For
$\xi=ga\cdot \xi_o\in \Xi_+$, with $g\in G$ and $a\in \exp (i\Omega_H)$,
define
\begin{equation}\label{real1} E_\R(\xi)=gNa\cdot x_o\subset
E(\xi )\, .
\end{equation}
Then $E_\R(\xi)$ is well defined, $G$-invariant
and a totally real submanifold of
$E(\xi)$. Further, the assignment $\Xi_+\ni \xi\mapsto E_\R(\xi)$
is  $G$-equivariant.

\par
Next we consider $S(z)\simeq K_\C/ M_\C$. Because of
the relation $K_\C =\zH^{-1} H_\C \zH$ there are  two natural real
forms. Accordingly we define for $z=ga\cdot x_o\in D$:
\begin{equation}\label{real2}
S^K_\R(z)=gaK\cdot \xi_o
\qquad \text{and}\qquad
S^H_\R (z)=ga \zH^{-1} H\zH \cdot \xi_o\, .
\end{equation}
Obviously $S^K_\R(z)$ and
$S^H_\R (z)$ are $G$-invariant totally real submanifold of $S(z)=gaK_\C\cdot x_o$
and the  maps $D\ni z\mapsto S^K_\R(z)$
and $D\ni z\mapsto S^H_\R(z)$  are
$G$-equivariant. Note that $S_\R^H(z)\simeq H/M$ as manifolds.

\subsection{Invariant measure on $Y$}\label{ss-measure}
Lemma \ref{nc} allows for a natural normalization
of the invariant measure on $Y$.
Assume that invariant measures on $G$, $A$ and $N$ have been fixed
and let $\cW_H=N_{K\cap H} (\fa)/ Z_{H\cap K}(\fa)$ be the
small Weyl group. By Lemma \ref{nc}
the union
$${\mathcal U}=\bigcup_{w\in \cW/\cW_H} AN w\cdot y_o$$
is disjoint  open and dense in $Y$. As the complement is an analytic set, it
follows that $Y\setminus {\mathcal U}$ has measure zero.
We  can normalize the invariant
measure
on $Y$ such that for all $f\in L^1 (Y)$:

$$\int_Y f (y) \, dy=\sum_{w\in \cW/ \cW_H} \int_A \int_N f
(anw\cdot y_o)\ da \ dn\, .$$

\section{The Frech\'et module $CR^\infty(\Xi_+)$}
\noindent In this section we use the right $A$-action on $\Xi_\R$
and $\Xi_+$ to define $G$-submodules of the smooth $A$-covariant
functions on $\Xi_\R$ respectively $CR$-functions on $\Xi_+$. Those
modules are the standard realization respectively a CR-realization
of the space of smooth vectors in the principal series
representations given by induction from the right. Note that later
we will use the induction from the left.

Recall, that $A$ acts on the space of horospheres from the right.
This action induces a right regular representation of $A$ on any
function space on $\Xi_\R$, $\Xi_+$ or any other right invariant set
of horospheres given by
$$(R(a)f)(\xi)=f(\xi\cdot a)\, .$$
Let $\rho=1/2 \sum_{\ga\in\Sigma^+} (\dim \fg^\ga)\cdot   \ga$. Let $\lambda
\in \fa_\C^*$. 

The
index ${}_\lambda$ will denote the subspace of
$(\lambda-\rho)$-covariant functions. In particular

$$C^\infty(\Xi_\R)_\lambda=\{ f\in C^\infty (\Xi_\R)\mid (\forall a\in A)
\ R(a)f =
 a^{\lambda-\rho} f\}\, .$$
We recall that $G$ acts $C^\infty(\Xi_\R)$ by left
translations in the argument
$$(L(g)f)(\xi)=f(g^{-1}\cdot \xi) $$
for $g\in G$, $f\in C^\infty(\Xi_\R)$  and $\xi\in \Xi_\R$. The so
obtained representation $(L, C^\infty(\Xi_\R)_\lambda)$ is the
smooth model of the spherical principal series with parameter
$\lambda$. 

\par Write $CR^\infty(\Xi_+)$ for the 
space of smooth $CR$-functions on $\Xi^+$
and set 
$$CR^\infty(\Xi_+)_\lambda=\{ f\in CR^\infty
(\Xi_+)\mid (\forall a\in A)\ R(a)f = a^{\lambda-\rho} f\}\, .$$ 
As characters on $A$
extend to holomorphic functions on $T(\Omega_H)$, it follows that
the restriction map

\begin{equation} \label{rest}{\rm Res}_\lambda:
CR^\infty(\Xi_+)_\lambda \to C^\infty(\Xi_\R)_\lambda, \ \
f\mapsto f|_{\Xi_\R}\end{equation}
is a $G$-equivariant topological isomorphism of $G$-modules.

\subsection{$CR$-realization of the $H$-spherical holomorphic vector}
For each $\lambda\in\fa_\C^*$ we define a certain $H$-invariant
element $f_\lambda\in CR^{-\infty}(\Xi_+)_\lambda$ which was called
the $H$-spherical holomorphic distribution vector in \cite{GKO}. The
generalized function $f_\lambda$ is defined by

$$f_\lambda(\xi) = a(\xi^{-1}z_H^{-1})^{\rho-\lambda} \qquad (\xi\in \Xi_+)\, .$$
We notice that on the dense subset

$$\Xi_+'=\bigcup_{w\in \cW} HwT(\Omega_H)\cdot\xi_o$$
of $\Xi_+$, the function belongs to $CR^\infty(\Xi_+')_\lambda$ and is
given by
$$f_\lambda(hwa\cdot\xi_o) = (w^{-1} z_H w)^{\lambda-\rho} a^{\lambda-\rho} \, .$$
For $\Re \lambda \ll 0$, this function is actually continuous on
$\Xi_+$ and the meromorphic continuation in $\lambda$ as a
distribution is achieved with Bernstein's theorem \cite{B}. There
are no singularities on the imaginary axis $i\fa^*$ 
and  we arrive
at a well defined analytic assignment

$$i\fa^*\to CR^{-\infty}(\Xi_+)_\lambda^H,  \ \lambda\mapsto f_\lambda\,$$
(cf. \cite{GKO}, Th. 2.4.1).

\section{The holomorphic horospherical Radon transform}\label{s=two}
\noindent
The real Radon transform on $X$ is the $G$-equivariant
injective map
\begin{equation}\label{rera}
\cR_\R: \cS(X)\to C^\infty(\Xi_\R), \ \ \cR_\R(f)(g\cdot \xi_o)=
\int_N f(gn\cdot x_o)\ dn  (\quad (g\in G))\, .\end{equation}
The purpose of this section is to show  that $\cR_\R$ has a natural
extension to a $G$-equivariant map

$$\cR : \cH^2(D)_0\to CR^\infty(\Xi_+)$$
which we call the holomorphic horospherical Radon transform.
Here $\cH^2(D)_0\hookrightarrow L^2 (X)$ is a
dense a subspace of the most-continuous Hardy space $\cH^2(D)\subset \cO(D)$ of $L^2(Y)$
(cf.\ \cite{GKO}).
On the infinitesimal level this extension is related to
the previously established fact (\ref{rest}), i.e. $C^\infty(\Xi_\R)_\lambda$ is canonically
$G$-isomorphic to $CR^\infty(\Xi_+)_\lambda$ via restriction.

\par This section is organized as follows: First we have to recall some
facts about the Fourier analysis on $X$, in particular Arthur's
spectral characterization of the Schwartz space $\cC(X)$.
Subsequently we give a brief summary on the most-continuous
Hardy space $\cH^2(D)$ of \cite{GKO}. Finally we define the holomorphic
horospherical Radon transform $\cR$ and discuss some of its properties.

\subsection{Fourier analysis on $X$}
We recall the compact realization
of the principal series representations.
Let $B=M\bs K$. For
$\lambda\in \fa^*_\C$  define a
representation $\pi_\lambda$ of $G$ on $L^2(B)$
by
\begin{equation}\label{eq-pi}
\pi_\lambda (g)f(Mk)=a(kg)^{\rho -\lambda}f(Mk(\emph{}kg))\, .
\end{equation}
Then $(\pi_\lambda ,L^2(B ))$
is unitary if $\lambda\in i\fa^*$.
We  write $\cH_\lambda =L^2(B)$ to
indicate the dependence of the $G$-action on $\lambda$.
Let $\fa^+=\{H\in \fa\mid (\forall \alpha \in \Delta^+ )\, \alpha (H)>0\}$ and
$$\fa_+^*=\{\lambda \in \fa \mid (\forall H\in\fa^+)\,  \lambda (H)>0\}\, .$$
Denote by $\hat{G}_{\mathrm{r}}$ the
reduced dual of $G$ and by $\hat{G}_{\mathrm{rsp}}$ the spherical
reduced dual. Then $i\fa_+^*\ni \lambda\mapsto
[\pi_\lambda ]\in \hat{G}_{\mathrm{rsp}}$
is an isomorphism of measure spaces. Here $[\pi_\lambda ]$ denotes the
equivalence class of $\pi_\lambda$. We have
\begin{equation}\label{eq-Fourier}
L^2(X )\simeq \int^\oplus_{i\fa^*_+}\cH_\lambda \,
\frac{d\lambda }{|\bc (\lambda )|^2}
\end{equation}
where $\bc (\lambda)$ is the Harish-Chandra $c$-function.
To explain the above isomorphism,
we need some basic fact on the Fourier transform
on $X$. For that, recall that $a:X\to A$ denotes  the $A$-projection
with regard to the horospherical coordinates 
$X=NA\cdot x_o\simeq N\times A$. Set
${\mathcal X}=B\times i\fa_+^*$ and define
a Radon measure $d\mu_{{\mathcal X}}$ on ${\mathcal X}$ by
$$d\mu_{{\mathcal X}}(b,\lambda ):=db
\frac{d\lambda}{|\bc(\lambda)|^2}\, .$$
For $f\in L^1(X)\cap L^2(X)$ define its
\textit{spherical Fourier transform} $\hat{f}:{\mathcal X}\to \C$ by
$$\hat{f} (b,\lambda)=\int_X f(x) a(bx)^{\rho-\lambda} \ dx\, .$$
We will also write $\cF_X(f)$ for $\hat{f}$.
We can normalize the left-Haar measure $dx$ on $X$ such that
the Fourier transform extends to an unitary isomorphism
$\hat{\hbox to 1em{}} : L^2(X)\to L^2({\mathcal X}, d\mu_{{\mathcal X}})$.
If $f$ is rapidly decreasing (see exact definition in a moment), then
the Fourier inversion formula holds pointwise:
$$f(x)=\int_{{\mathcal X}} \hat f(b,\lambda) a(bx)^{\rho+\lambda}
\ d\mu_{{\mathcal X}}(b,\lambda )\qquad (x\in X)\, .$$ For $\lambda
\in  i\fa^*$ define $\hat{f}_\lambda \in L^2(B)$ by $b\mapsto
\hat{f}_\lambda (b)=\hat{f}(b,\lambda)$. Then the isomorphism in
(\ref{eq-Fourier}) is given by
$$L^2(X)\ni f\mapsto (\hat{f}_\lambda )_\lambda \in \int^\oplus_{i\fa^*_+}\cH_\lambda \,
\frac{d\lambda }{|\bc (\lambda )|^2}\, .$$

In the following we will also need the operator valued Fourier transform.
If $\cH$ is a Hilbert space, then
$B_2(\cH)\simeq \cH \hat \otimes \cH^*$ denotes
the Hilbert space of Hilbert-Schmidt operators
on $\cH$.
Write
$$L^2(G)_{\rm sph}=\int_{i\fa_+^*}^\oplus B_2(\cH_\lambda) \
\frac{d\lambda}{|\bc(\lambda)|^2}$$
for the $K$-spherical
spectrum in $L^2(G)$.
Recall that the isomorphism is given by the operator valued
Fourier transform $\cF (f)(\lambda )= \int_G f(x)\pi_\lambda (x)\, dx$,
$f\in L^1(G)\cap L^2(G)$. The inverse map is
$$f(g)=\int_{i\fa_+^*}\Tr (\pi_\lambda (g^{-1})\cF (f)(\lambda ))\,
\frac{d\lambda }{|\bc (\lambda )|^2}\, .$$

The constant function $v_{K,\lambda}=\one_{B}$ defines a normalized $K$-fixed
vector in $\cH_\lambda$.
Assume that $f\in L^1(G)\cap L^2(G)_{\rm sph}$. Then, because
$\cF (f)(\lambda)=\cF (R_k f)(\lambda)=\cF (f)(\lambda)\pi_\lambda (k)$, it follows, that
\begin{equation}\label{eq-FT}
\cF (f)(\lambda)v=
\la v,v_{K,\lambda}\ra \cF(f)(\lambda)v_{K,\lambda}=
\la v,v_{K,\lambda} \ra \hat{f}_\lambda \, .
\end{equation}

For $x=k_1\exp (Z) k_2\in G$, with $Z\in \fa$ and $k_1,k_2\in K$, let
$\sigma (x)=-B (Z,d\theta(\bf 1) (Z))$, where $B$ is the Killing form on $\fg$.
Denote by $U(\fg )$ the Universal enveloping algebra of $\fg$ and by
$\varphi_0$ the basic spherical function.
For $D,E\in U(\fg)$,  $s\in \R$ and $f\in C^\infty (G)$, let
$$p_{D,E,s}(f):=\sup_{x\in G}|L_D R_Ef(x)|\varphi_0(x)^{-1}(1+\sigma (x))^s\, .$$
Then $\cC (G)$ is the space of smooth functions on $G$ such
that $p_{D,E,s}(f)<\infty $ for all such $E,D$ and $s$, cf. \cite{HC66}, \S 9
.
We set
$$\cC (G)_{\rm sph}= \cC (G)\cap L^2(G)_{\rm sph}\, .$$

For $\tau\in \hat K$ we write $|\tau|$ for the norm of the corresponding 
highest weight. If $\tau\in \hat K$, then we write $L^2(B)_\tau$ for the subspace of $L^2(B)$ which transform under 
$\tau$. We denote by $\D(i\fa^*)$ the algebra of all 
constant coefficient differential operators on $\fa^*$
and set $\cS(i\fa_+^*)=\{f|_{i\fa_+^*}\mid
f\in \cS (i\fa^*)\}$.

We cite a theorem of Arthur \cite{Arthur}, p. 4719, specialized to the spherical case:
\begin{theorem}\label{th-arthur}
The operator valued Fourier transform $\cF$ 
is a topological linear isomorphism  from 
$\cC(G)_{\rm sph}$ onto
\begin{align*} \{ &A(\, \cdot \, )\in \int_{i\fa_+^*}^\oplus B_2(\cH_\lambda) \,
\frac{d\lambda}{|\bc(\lambda)|^2} \mid (\forall v,w\in L^2(B )_{\rm K-fin})
\, \la A(\cdot )v,w\ra
\in \cS(i\fa_+^*)\, , \\ 
&\forall D\in \D(i\fa^*),  \forall n\in \N_0 
 \sup_{\lambda\in i\fa_+^*, \sigma,\tau \in \hat K\atop v\in L^2(B)_\sigma,  w\in L^2(B)_\tau}
{|D \la A(\lambda) v, w\ra|\over \|v\|\cdot \|w\|} (1+|\lambda|)^n (1+ |\tau|)^n 
(1+|\sigma|)^n<\infty \}\, .\end{align*}
\end{theorem}

\subsection{The most-continuous Hardy space}
We recall now the spectral definition of
the Hardy space $\cH^2(D)$ from \cite{GKO}.
For $v\in \cH_\lambda$
define an analytic function
$f_{v, \lambda}$ on $X$ by

$$f_{v,\lambda}(x)=\langle \pi_\lambda(x^{-1}) v, v_{K, \lambda}\rangle=
\langle  v, \pi_\lambda (x) v_{K, \lambda}\rangle\, .$$
Let us denote by $z\mapsto \bar{z}$ the complex conjugation 
in $G_\C$ with respect to $G$
and  recall that $f_{v,\lambda}$ extends to a holomorphic
function $\tilde f_{v,\lambda} $ on $D$ via
$$\tilde f_{v,\lambda}(x)=\langle v, \pi_\lambda(\bar{x})v_{K,\lambda}
\rangle$$
for $x\in D$, cf. \cite{GKO} and \cite{GKO2}, Proposition 2.2.3. 
In particular
$$f_{v,\lambda}(ga \cdot x_o)=  \langle \pi_\lambda(g^{-1}) v, \pi_\lambda (a^{-1})
v_{K,\lambda }\rangle $$
for $g\in G$ and $a\in \exp (i\Omega_H)$

Define a generalized hyperbolic
cosine function on $i\fa^*$ by
\begin{equation}\label{eq-cosh}
\COSH (\lambda)=\sum_{w\in \cW/\cW_H} \zH^{2w^{-1}\lambda}
\end{equation}
for $\lambda\in i\fa^*$.
Define a measure $\mu$ on $i\fa_+^*$ by
\begin{equation}\label{eq-mu}
d\mu(\lambda) =\frac{d\lambda} {\COSH(\lambda) \cdot |\bc(\lambda)|^2}\, .
\end{equation}
With this preparation we can define the unitary representation $(L ,\cH^2(D))$
of $G$ by
$$(L ,\cH^2(D))=\int^\oplus_{i\fa_+^*}(\pi_\lambda ,\cH_\lambda )\,
d\mu (\lambda)\, .$$
Thus $\cH^2(D)$ is the Hilbert space of all measurable functions
$s :i\fa_+^*\to L^2(M\bs K)$ such that
$\|s\|^2=\int_{i\fa_+^*}\|s (\lambda ) \|^2 \, d\mu (\lambda )<\infty$.
In the sequel we often write  write $s_\lambda$ for
$s (\lambda )$. Let us denote by $\|\cdot\|_H$ the norm on $\cH^2(D)$.
Recall from \cite{GKO} that the map
$$\Phi: \cH^2(D)\hookrightarrow \cO(D), \ \ s=(s_\lambda)\mapsto \left(
x\mapsto \int_{i\fa_+^*} \tilde f_{s_\lambda, \lambda}(x) \, d\mu(\lambda)\right)$$
is a $G$-equivariant continuous injection. In the sequel we often view
$\cH^2(D)$ as a subspace of $\cO(D)$; we 
call
$\cH^2(D)$ the {\it most-continuous Hardy space of $Y$}.
This notion is motivated by the  main result of \cite{GKO} which states that there exists
an isometric $G$-equivariant value mapping
$$b: \cH^2(D)\to L_{\rm mc}^2 (Y), \ \ f\mapsto b(f)\, $$
which is isometric onto a multiplicity free subspace 
of $L_{\rm mc}^2(Y)$. 

\subsection{The Fourier Transform on $X$ and the Hardy space}
The definition of $\cH^2 (D)$ in the previous subsection does not
use the Fourier transform on $X$.  But the following Lemma
shows that the space $\cH^2 (D)$ has a natural description in
terms of the Fourier transform.

\begin{lemma}\label{eq=hn} Let $f\in \cH^2(D)$. Then the following 
assertions hold:
\begin{enumerate}
\item $f|_X\in L^2(X)$ and
\begin{eqnarray*}
f (z) &=&\int_{i\fa_+^*} \widehat{f|_X}(b,\lambda ) \COSH (\lambda)
\, a(bz)^{\lambda +\rho}\, d\mu (b,\lambda )\qquad (z\in D)\\
\|f\|_H^2&=&\int_{\mathcal X} |\widehat{f|_X}(b,\lambda)|^2\cdot
\COSH(\lambda)\, d\mu_{\mathcal X}(b,\lambda)\ge \|f|_X\|_{L^2(X)}^2
\, .
\end{eqnarray*}
\item If $f=\Phi^{-1}(s_\lambda )$, then
$$(\widehat{f|_X})_\lambda = \frac{s_\lambda}{\COSH (\lambda )}\, .$$
\item For $a=\exp (iY)\in \exp (i\gO_H)$ let  $f_a:G\to \C$ by
$f_a(g)=f(ga\cdot x_o)$. Let $Q\subset \exp (i\gO_H)$ be compact
and such that $a\in Q$. Then there exists a constant $C_Q>0$ such that
$$ \|f_a\|_{L^2(G)}\le C_Q\| f\|_H $$
\end{enumerate}
\end{lemma}

\begin{proof} (1) and (2).  Let $f\in \cH^2(D)$ and
$f=\int_{i\fa_+^*} \tilde f_{s_\lambda,\lambda} \, d\mu (\lambda )$.
Then obviously we have (2), i.e.,
\begin{equation}\label{eq-decomposition}
f=\int_{i\fa_+^*}\frac{\tilde f_{s_\lambda,\lambda} }{\COSH (\lambda )}\,
\frac{d\lambda}{|\bc (\lambda )|^2}
\end{equation}
and
\begin{eqnarray*}
\|f\|_H^2&=&\int_{i\fa_+^*}\|s_\lambda \|^2_{L^2(B)}\, d\mu (\lambda )\\
&=& \int_{i\fa_+^*}\left\|\frac{s_\lambda}{\COSH (\lambda )} \right\|^2_{L^2(B)}
\, \COSH (\lambda )\, \frac{d\lambda }{|\bc (\lambda )|^2}\\
&\ge &\int_{i\fa_+^*}\left\|\frac{s_\lambda}{\COSH (\lambda )} \right\|^2_{L^2(B)}
\, \frac{d\lambda }{|\bc (\lambda )|^2}\\
&=& \|f|_X\|^2_{L^2 (X)}\, .
\end{eqnarray*}
Thus $f|_X\in L^2(X)$ and we can write
$f|_X=\int_{\mathcal X} \widehat{f|_X}(b,\lambda)
a(b\, \cdot
\, )^{\rho+\lambda}
db \frac{d\lambda}{ |\bc(\lambda)|^2}$.
Equation (\ref{eq-decomposition}) implies that
$\widehat{f|_X} (b,\lambda )=s_\lambda (b)/\COSH (\lambda )$
or $s_\lambda (b)=\widehat{f|_X}(b,\lambda )\COSH (\lambda )$
for almost all $\lambda$. This finishes
the proof of (1) and (2).

\par
(3) We recall  Faraut's version of
the  Gutzmer identity \cite{Faraut} \begin{equation}\label{eq=Gutzmer}
\int_G |f(ga\cdot x_o)|^2 \ dg =\int_{\mathcal X} |\widehat{ f|_X}
(b,\lambda)|^2
\varphi_\lambda(a^2) \ db \frac{d\lambda}{|\bc(\lambda)|^2}
\end{equation}
where  $\varphi_\lambda(a^2)$
is the analytically continued spherical function given by
$$\varphi_\lambda(a^2)=\int_K \left|a(ka)^{\rho+\lambda}\right|^2 \ dk \,$$
(cf. \cite{KS}, Sect. 4).
Now for a compact subset $Q\subset \exp (i\gO_H)$ there exists
a constant $C_Q>0$ such that
\begin{equation}\label{EST} (\forall \lambda\in i\fa^*)\qquad 
\varphi_\lambda(a^2)\leq C_Q\COSH(\lambda)
\end{equation}
for all $\lambda$   (cf.\ \cite{kos05}, Lemma 2.1)
and the assertion of the lemma follows.
\end{proof}

In order to define the Radon transform for functions in the Hardy space
we first need a technical fact, interesting on its own.

Let $
\left(\int^\oplus_{i\fa^*_+}\cH_\lambda \, d\mu (\lambda )\right)_0$
denote the space of all  sections $(s_\lambda)$ such that
for all $v\in L^2(B)_{\hbox{$K$-finite}} $
$$i\fa_+^*\ni \lambda \mapsto
\la s_\lambda ,v \ra \in \cS(i\fa_+^*)\, $$
and 
\begin{equation}\label{rd}
(\forall D\in \D(i\fa^*), \forall n\in \N_0) \sup_{\lambda\in i\fa_+^*, \tau\in \hat K \atop v\in L^2(B)_\tau}
{| D \la s_\lambda, v\ra |\over \|v\|} (1+|\lambda|)^n (1+|\tau|)^n<\infty
\, .\end{equation}

Then we set
\begin{equation}\label{eq-Hc0}\cH^2(D)_0=\Phi^{-1}
\left( \left(\int^\oplus_{i\fa^*_+}\cH_\lambda \, d\mu (\lambda )\right)_0 \right)\, .
\end{equation}

\begin{theorem}\label{th=Schwartz} Let $f\in \cH^2(D)_0$. Fix
$z\in T(\Omega_H)$. Then the function
$$G\ni g\mapsto f(gz\cdot x_o)\in\C$$
belongs to $\cC (G)$. Moreover, the following functions
are locally bounded on $T(\Omega_H)$:
\begin{enumerate}
\item $z\to \int_G |f(gz\cdot x_o)|^2 \ dg$
\item  $z\to \int_N |f(nz\cdot x_o)| \ dn$
\end{enumerate}
\end{theorem}

\begin{proof}
Without loss of generality we may assume that
$z=a\in \exp(i\Omega_H)$. In the sequel we often identify $f$ with $f|_X$, 
a right $K$-invariant function on $G$.
Let $v,w\in L^2(B)_{K-{\rm finite}}$. Then, 
by (\ref{eq-FT}) and Lemma \ref{eq=hn}:
\begin{eqnarray*}
\la  \cF (R_a f)(\lambda )v,w \ra & =&\la  \cF (f)(\lambda )\pi_\lambda (a)v, w\ra\\
&=&\la \pi_\lambda (a)v,v_{K,\lambda } \ra \la \hat{f}_\lambda ,w\ra\\
&=&\frac{\la v, \pi_\lambda (a)v_{K,\lambda } \ra}{\COSH (\lambda )}\la s_\lambda ,w\ra \, .
\end{eqnarray*}

Let us write $F(\lambda, v,w):=\la  \cF (R_a f)(\lambda )v,w \ra $. It is 
clear that $F$ is smooth in the $\lambda$-variable.  In order to show that 
$g\mapsto f(ga\cdot x_0)$ belongs to $\cC(G)$ we use 
Arthur's Theorem \ref{th=Schwartz}: 
we have to show for all $n\in \N_0$ and $D\in \D(i\fa^*)$ that

\begin{equation}\label{p1}\sup_{\lambda\in i\fa_+^*, \sigma,\tau\in \hat K\atop 
v\in L^2(B)_\tau, w\in L^2(B)_\sigma} {\left|D F(\lambda, v, w) \right|\over \|v\|\cdot\|w\|}
(1+|\lambda|)^n (1+|\tau|)^n (1+|\sigma|)^n <\infty\, .\end{equation}

We use the expression for $F(\lambda, v,w)$ derived above to 
compute its derivatives: For fixed $D$, we find according to Leibniz
a finite set of $D_1, D_1'\ldots, D_N, D_N'$ in $\D(i\fa^*)$ such that

\begin{equation}\label{p2}
D F(\lambda, v, w)= \sum_{j=1}^N 
\left[D_j 
\frac{\la v,\pi_\lambda (a)v_{K,\lambda } \ra}{\COSH (\lambda )}\right]
\cdot \left[D_j' \la s_\lambda, w\ra\right]\, .
\end{equation}

According to (\ref{rd}), for all $D\in \D(i\fa^*)$, 
$n\in \N_0$ the following estimate 
holds: 
\begin{equation}\label{p3}\sup_{\lambda\in i\fa_+^*, \sigma\in \hat K\atop 
w\in L^2(B)_\sigma} {\left|D \la s_\lambda, w\ra\right |\over \|w\|}
(1+|\lambda|)^n (1+|\sigma|)^n <\infty\, .\end{equation}
Now, after combining (\ref{p1})-(\ref{p3}) we are left to show 
for all $D\in \D(i\fa^*)$, $n\in\N_0$: 

 \begin{equation}\label{p4}\sup_{\lambda\in i\fa_+^*, \tau\in \hat K\atop 
v\in L^2(B)_\tau} {\left|D \frac{\la v, \pi_\lambda (a)v_{K,\lambda } \ra}{\COSH (\lambda )}
\right|\over \|v\|}
(1+|\tau|)^n  <\infty\, .\end{equation}

To establish (\ref{p4}), we first note that 
$\pi_\lambda(a)v_{K,\lambda }$ is an analytic vector 
for the representation. 

\par Now if $w\in \cH_\lambda^\omega= C^\omega(M\bs K)$ 
and $w=\sum_{\tau \in \hat K} w_\tau$ is its expansion in $K$-types, 
then we recall from  \cite{KV}, Th. 2.2 (1), that there 
exists a $\delta>0$ such that 
$$\|w_\tau\| << e^{-\delta|\tau|}\, .$$
To be more precise, given $\delta>0$ sufficiently small, 
there exists a ball $U$  
arround ${\bf 1}$ in $K_\C$ such that $\pi_\lambda(t)w$ 
exists for all $t\in U$ and 
$$\|w_\tau\| 
\leq (\sup_{t\in U} \|\pi_\lambda (t) w\|)
\cdot  e^{-\delta/2 |\tau|}\, .$$

\par Coming back to our original situation where our 
analytic vector is $w=\pi_\lambda(a) v_{K_\lambda}$,  
we obtain that $\pi_\lambda(t)\pi_\lambda(a) v_{K_\lambda}$ 
exists for all $t\in U$.
It follows that $\pi(G U aK_\C) v_{K_\lambda}$ exists. As the crown is 
a domain of holomorphy for principal series 
representations (cf. \cite{K}), we obtain,  by continuity,  
a compact neighborhood $Q$ of $a$ 
in $\exp(i\Omega_H)$ such that $UaK_\C \subset G Q K_\C$. 
 By the very nature of 
the principal series, it is easy to see that 
$U$ can be chosen independent from $\lambda$.
Hence we get a $C>0$ such that for all $\lambda\in i\fa^*$, $\tau\in
\hat K$ and $v\in L^2(B)_\tau$
\begin{equation}\label{p5}
{|\la v, \pi_\lambda (a)v_{K,\lambda }\ra| 
\over \|v\|} \leq C (1+|\tau|)^{-n}\cdot \sup_{b\in Q} 
\|\pi_\lambda (b)v_{K,\lambda }\|\, .
\end{equation}
Moreover, 
$\sup_{b\in Q}{\|\pi_\lambda (b)v_{K,\lambda }\|\over 
\sqrt{\COSH (\lambda )}}$ is uniformly bounded 
in $\lambda$ by (\ref{EST}). This, in combination with (\ref{p5}), 
proves (\ref{p4}) for the case of 
$D={\bf 1}$. Now,  if $D$ is  of order 
$k$,  we first 
observe that 

$$(\forall k\in K)(\forall b\in T(\Omega_H))\quad 
\left[D \pi_\lambda (b)v_{K,\lambda }\right](Mk)=
p_k ((\lambda +\rho)(\log a(kb))   \left(\pi_\lambda (b)v_{K,\lambda }\right)(Mk)$$
for a polynomial $p_k$ of degree $k$ and independent from $\lambda$. 
It follows from the complex convexity theorem (\ref{cc})  that 
$$\sup_{b\in Q} \sup_{k\in K} |p_k (\lambda +\rho)(\log a(kb))| \leq C \cdot 
\sqrt{\COSH (\lambda )}$$
Furthermore the decay of ${1\over \COSH (\lambda )}$ is not affected 
by differentiation. Therefore we obtain (\ref{p4}) for all $D$. 

\par Moving on to (1) and (2), we observe that 
statement (1) is Lemma \ref{eq=hn}, part 3;  part (2) follows  from 
the already established fact in conjunction with Lemma 22 in \cite{HC66}.
\end{proof}

\subsection{The definition of the Radon Transform}
Denote by $CR(\Xi_+)$ the vector space of continuous $CR$-functions
on $\Xi_+$, i.e. the space of continuous  functions
on $\Xi_+\simeq K/M\times T(\Omega_H)$ (cf. (\ref{eq=decomp}))
which are holomorphic in the second variable.

\begin{lemma}\label{prele} Let $f\in \cH^2(D)_0$. Then the assignment
$$\Xi_+\ni \xi=ga\cdot\xi_o\mapsto  a^{-2\rho}\int_N f(gna\cdot x_o)\ dn\in \C
\qquad (g\in G, a\in \exp(i\Omega_H))$$
defines a $CR$-function on $\Xi_+$.
\end{lemma}

\begin{proof} It follows from Theorem \ref{th=Schwartz} that the right hand
side is a continuous function. It remains to show that is a $CR$-function.
For that let $g=kb$ for $k\in K$ and $b\in A$. The right hand
side becomes
$$a^{-2\rho}\int_N f(kbna\cdot x_o)\ dn=(ab)^{-2\rho}\int_N f(knba\cdot x_o)\ dn$$
and the holomorphicity in $ab$  follows with Theorem \ref{th=Schwartz}.
\end{proof}
In view of this lemma, the prescription

$$\cR: \cH^2(D)_0\to CR(\Xi_+), \ \ f\mapsto
\left( \xi=ga\cdot\xi_o\mapsto a^{-2\rho}\int_N f(gna\cdot x_o)\ dn\right)$$
is a well defined and continuous $G$-equivariant map.
We call $\cR$ the \textit{ holomorphic horospherical Radon transform}.

\begin{remark} (a) We recall the real Radon transform on $X$ from (\ref{rera}).Now if $f\in \cH^2(D)_0$, then $f|_X\in L^2(X)$ by Lemma \ref{prele}. 
Thus, for $g\in G$, we get 
$$\cR(f)(g\cdot \xi_o)=\cR_\R(f|_X)(g\cdot \xi_o)\, . $$
In other words, the holomorphic Radon transform restricted 
to $\Xi_\R$ agrees with the real 
Radon transform of the restricted function on $X$.
\par (b) Notice that $E(\xi)\cap D$ for $\xi\in \Xi_+$
contains the real horosphere $E_\R(\xi)$. The
holomorphic Radon transform $\cR$ then writes as
$$\cR(f)(\xi)=\int_{E_\R(\xi)} f(\xi') \ d\nu_\xi(\xi')$$
with $d\nu_\xi$ equals $a^{-2\rho}$-times the measure on $E_\R(\xi)$ 
obtained by
the natural identification
of the real horosphere $E_\R(\xi)$ with
$N$. It is clear that any other $N$-orbit in $E(\xi)\cap D$
would yield the same result.

\end{remark}

\par
If the function $f\in \cH^2(D)_0$ is left
$K$-invariant, then we  can define the \textit{holomorphic Abel-transform}
by:

$$\cA(f)(z)=z^{-\rho} \int_N f(nz\cdot x_o)\ dn \qquad (z\in T(\Omega_H) )\, .$$
Note that $\cA$ is just the restriction of the holomorphic
Radon transform to $K$-invariant function
(modulo the $z^{-\rho}$-factor).
Further let us remark that $\cA$ gives a continuous mapping

$$\cA:\cH^2(D)_0^K\to \cO(T(\Omega_H))^\cW, \ \
\cA(f)(z)=z^{-\rho} \int_N f(nz\cdot x_o)\ dn\ .$$

\section{The holomorphic Radon transform as Cauchy integral I: The hyperboloid}
\noindent In this and the next section we will show (for an
appropriate class of $Y$'s) that the holomorphic Radon transform
on the NCC-space $Y$ can
be expressed as a Cauchy type integral. For that purpose
it is instructive to explaining the example of the hyperboloid first. For
earlier treatments of the hyperboloid with alternative methods we
refer to \cite{G1}, \cite{G2}. We start by recalling some standard
function spaces on $Y$.

\subsection{Function spaces} Let $Y$ be a NCC-space. For $g\in G$ let
$\Theta (g)=\varphi_0(g\tau (g)^{-1})^{1/2}$. Then $\Theta$ is left $K$-invariant
and right $H$-invariant. For $g=k\exp (Z) h$ with
$Z\in \fa$ define
$\|g\cdot y_0\|:=\|Z\|$, where $\|Z\|=\sqrt{\Tr \ad (Z)^2}$.
Let $D\in U(\fg )$, were $U(\fg )$ is the enveloping algebra
of $\fg_\C$, and view $L_D$ as a differential operator on
$Y$. For $n\in \N$ and $f\in C^\infty (Y)$ define
$$p_{n,D}(f)=\sup_{y\in Y}\Theta_G(y)^{-1}(1+\|y \|)^n |L_Df(y)|\, .$$
Then the Schwartz space $\cC (Y)$ is defined as the space
of smooth function on $Y$ such that $p_{n,D}(f)<\infty $ for
all $n$ and $D$. It is well known, that $\cC (Y)\subset L^2(Y)$,
but $\cC (Y)$ is not contained in $L^1(Y)$. We will therefore
need a smaller space to make sure that the Cauchy integral exists.
For that, define for $r>0$ the space
$$\cS_r(Y):=\{ f\in C^\infty (Y)\mid (\forall D\in U(\fg ))\, \sup_{y\in Y}e^{r\|y\|}|L_Df(y)|
<\infty\}$$
and
$$\cS (Y):=\bigcap_{r>0}\cS_r(Y)\, .$$
Then $\cS (Y)\subset L^1(Y)\cap L^2(Y)$ and
$C_c^\infty (Y)\subset \cS (Y)\subset \cC (Y)$.
The space $\cS (Y)$ is called the (zero) \textit{Schwartz space},
cf. \cite{FD91}. It follows from Theorem 3 in \cite{FD91}
and our spectral definition of the space $\cH^2(D)$ that
$$\cH^2 (D)_{00}:=\{ f\in \cH^2(D)_0\mid 
f\ \hbox{$K$-finite}, \ b(f)\in \cS(Y)$$
is dense in $\cH^2(D)$. Here 
used the fact 
that elements in $\cH^2(D)_0$ have
boundary values on $Y$ (cf.\ \cite{GKO}, Sect. 3 for a similar argument). In
particular we have that every element $f\in \cH^2(D)_{00}$ is
integrable on $Y$, has a holomorphic extension to $D$, and
that the integral
$\int_N f(anw\cdot y_0)\, dn$ is well defined for all
$a\in A$ and $w\in \cW$. We will use this without comments
in the following.

\subsection{The Radon transform and Cauchy integral on the hyperboloid}
Assume that
$n\ge 2$ and let $G=\SO_e(1,n)$ be the Lorentz group.
Let us fix our choices
for the groups $A$, $N$ and $K$.
For the  maximal
compact subgroup we take
$$K=\left\{k_R=\begin{pmatrix} 1 & 0\\
0 &R\end{pmatrix}\mid R\in \SO(n)\right\}\simeq \SO (n)\, .$$
Next, for $z\in \C$ we set

$$a_z=\begin{pmatrix}  \cosh z & 0 & \sinh z\\
0 & \one_{n-1}& 0\\
\sinh z  & 0 & \cosh z \end{pmatrix}$$
and
$$A=\{ a_t\mid t\in \R\}\qquad \hbox{and}\qquad A_\C=\{
a_z\mid z\in \C\}\, .$$
Note that $\fa=\R Z$ and $a_z=\exp (z Z )$ with $Z=E_{1\, n+1}+E_{n+1\, 1}$.
The only positive root is
$\alpha$, determined by $\alpha (Z)=1$ and hence  $\ZH=\frac{\pi }{2}Z$.
We also have, that $\rho = \frac{n-1}{2} \alpha $.

Further, for $v\in \C^{n-1}$ and $(v,  v) =\sum_{j=1} v_j v_j$ we define
an unipotent matrix
$$n_v=\begin{pmatrix}  1+ \frac{1}{2} (v, v) & v^T & -\frac{1}{2} (v, v)\\
v & \one_{n-1}& -v \\
\frac{1}{2} (v, v) & v^T & 1 -\frac{1}{2} (v,  v)
\end{pmatrix}\, . $$
Then $N$ and $N_\C$ are given by
$$N=\{ n_v\mid v\in \R^{n-1}\}\qquad \hbox{and}\qquad N_\C=\{
n_v\mid v\in \C^{n-1}\}\, .$$

\par
Define a quadratic form
$$z,w\mapsto z\cdot w=z_0w_0-\sum_{j=1}^nz_jw_j$$
on $\C^{n+1}$ and let
$$\square(\z)=z_0^2-z_1^2-\ldots -z_n^2 \qquad (\z=(z_0, \ldots, z_n)^T\in \C^{n+1})$$
be the corresponding quadratic form.
We define the real and complex hyperboloids
by
$$X=\{ \x\in\R^{n+1}\mid \square(\x)=1, x_0>0\}$$
and
$$X_\C=\{ \z\in\C^{n+1}\mid \square(\z)=1\}\, . $$
As a common base point for $X$ and $X_\C$ we take
$\x_o=(1,0,\ldots, 0)^T$ and note that
the map
$$G_\C/ K_\C \to X_\C, \ \ gK_\C\mapsto g(\x_o)$$
is a $G_\C$-isomorphism which identifies $G/K$ with $X$.

We have that $a_z\cdot \x_o=(\cosh (z),0,\ldots ,0,\sinh(z))^T$
and hence $\y_0=\zH\cdot \x_o=(0,\ldots ,0,i)^T$. It is clear that the stabilizer
of $\y_o$ is
$$H=\left\{\begin{pmatrix} h & 0\\
0 & 1\end{pmatrix}\mid h\in \SO_e(1,n-1)\right\}\simeq \SO_e(1,n-1)\, .$$
We have therefore with this identification:
$$D=\{ \z=\x+i\y\in X_\C\mid \square(\x)>0, x_0>0\}$$
and
$$Y=G (\y_0 )=\{ i\y\in i\R^{n+1}\mid \square(\y)=-1\}\, . $$

Set
$$\Xi=\{ \zeta\in \C^{n+1}\mid \zeta\neq 0, \square(\zeta)=0\}\, .$$
If $\xi_o=(1,0, \ldots, 0, 1)^T\in \Xi$, then
the stabilizer of $\xi_o$ is $M_\C N_\C$ and  the map
$$G_\C/ M_\C N_\C \to \Xi, \ \ gM_\C N_\C\mapsto g\cdot\xi_o$$
is a $G_\C$-isomorphism. Now
the $CR$-submanifold $\Xi_+\subset\Xi$ is described as

\begin{eqnarray*} \Xi_+&=&G\left\{ (e^{it}, 0, \ldots, 0, e^{it})^T\mid |t|<\frac{\pi}{2}, t\in \R\right\}\\
&=&\{\zeta=\xi+i\eta\in\Xi: \square(\xi)=\square(\eta)=0; \ \xi\neq 0\}_0\, .
\end{eqnarray*}

We will also use certain $G$-subdomains of $\Xi_+$: for $0<c\leq \frac{\pi}{ 2}$ set
$$\Xi_c=G\left\{ (e^{it}, 0, \ldots, 0, e^{it})^T\mid |t|<c, t\in \R \right\}\, .$$

\par
In order to define the Cauchy-transform we need to establish
a simple, but important, technical fact.

\begin{lemma}\label{lem=hyp} For all $\xi\in \Xi_+$ and $y\in Y$ one has
$$\xi\cdot y\not\in \R^\times \, .$$
More precisely, for all $0<c<\frac{\pi}{2}$ there exists a $C>0$ such that
$$(\forall y\in Y) (\forall \xi\in \Xi_c) \qquad |1-\xi\cdot y|>C\, .$$
\end{lemma}

\begin{proof} By $G$-equivariance of the form we may assume that
$\xi=e^{it}\xi_o$. Let $y=i\y$ for $\y\in \R^{n+1}$. Thus
$\xi\cdot y= ie^{it} \xi_o\cdot \y$ with
$\xi_o\cdot\y\in \R$. As $|t|<\frac{\pi }{2}$, the assertions follow.
\end{proof}

We now define the Cauchy-kernel function
$$\cK: \Xi_+\times Y\to \C , \ \ (\xi, y)\mapsto
\frac{1}{1- \xi\cdot y}\, .$$
In view of the previous lemma, this function is
defined, continuous and bounded on all subsets
$\Xi_c\times Y$ for $c<\frac{\pi}{2}$.  Moreover, $\cK$ is a $CR$-function in
the first variable and $G$-invariant, i.e., $\cK(g\cdot \xi ,g\cdot y)=\cK(\xi,y)$.
In particular, the function
$$G\ni g\mapsto \cK (g)=\cK (\xi_o,g\cdot \y_0)\in \C$$
is
left $N$-invariant and right $H$-invariant, a fact that we will use in a moment.
We will therefore identify $\cK$ with a function on $Y$ without
further comments.
A simple calculation shows that
\begin{equation}\label{eq-KSO}
\cK(g)=\frac{1}{1-i(g_{0\, n}-g_{n\, n})}\, .
\end{equation}

We have $\cW_H=\{\one \}$, and
$\cW=\{\one , \epsilon\}$ where $\epsilon=-1$ on $\fa$.
As $\epsilon$ corresponds to the matrix
$w=\begin{pmatrix} I_{n-1} & 0\cr 0&-I_2\end{pmatrix}$
it follows that
\begin{equation}\label{eq-K2}
\cK (a_z)=\frac{1}{1-ie^{-z}}\qquad\text{and}\qquad
\cK (a_zw)=\frac{1}{1+ie^{-z}}\, .
\end{equation}

Write $\cS(Y)$ for the Schwartz-space on $Y$.  Henceforth we will
make the assumption that $n$ is even and define the {\it Cauchy -
transform} by

$$\cC: \cS(Y)\to CR(\Xi_+), \ \ \cC(f)(\xi)=\int_Y
\frac{f(y)}{1- \xi\cdot y}\, dy= \int_Y f(y)K(\xi,y)\, dy\, $$
where $dy$ denotes the  $G$-invariant measure on $Y$ from Subsection 1.6. 

\begin{remark} In the following we have to raise complex powers
of  elements $z\in A_\C$ which are in the boundary of 
$\exp(i\Omega_H)$. This is not a problem  as the map 

$$\exp: 2\Omega_H\to A_\C , \ \ Y\mapsto \exp(iY)$$  
is injective; a fact which holds in full generality. 
\end{remark}

\begin{theorem} Let $G=\SO(n,1)$ with $n=2k$ even. Let $f\in \cH^2(D)_{00}$. 
Then, up to normalization of measures on $A$, 
one has 
$$\cC(f)(\xi)=(-1)^{k-1} 2\pi\cdot \cR(f)(\xi) \qquad (\xi\in \Xi_+)\, .$$
\end{theorem}

\begin{proof} Since both $\cC(f)$ and $\cR(f)$ are $CR$-functions, it
is sufficient to show that both coincide on $G/MN\subset \Xi$.
Moreover, by $G$-equivariance of both maps, it is in fact sufficient to show that
\begin{equation}\label {e1}  \cC(f)(\xi_o)=(-1)^{k-1}
2\pi \cdot \cR(f)(\xi_o)\, .\end{equation}
Using (\ref{eq-K2}) and that $\zH^{2\rho}=(-1)^{k-1} i$
and $z_H^{-2\rho}=\overline{\zH^{2\rho}}=(-1)^k i$ we get:

\begin{eqnarray*} \cC(f)(\xi_o) &=& \int_Y  f(y) \cdot \cK( y)\ dy \\
&=& \sum_{w\in \cW} \int_A \int_N f(anw\cdot \y_o)\cdot  \cK( anw)\ dn \ da \\
&=& \sum_{w\in \cW} \int_A \int_N f(anw\cdot \y_o) \cdot \cK( aw)\ dn \ da \\
&=& \sum_{w\in \cW} \int_A \int_N f(anw\zH \cdot \x_o) \cdot \cK(aw)\ dn \ da \\
&=& \sum_{w\in \cW} \int_A \int_N f(an\zH^w \cdot \x_o)\cdot  \cK( aw)\ dn \ da \\
&=& \sum_{w\in \cW} \int_A \int_N f(an\zH^w  \cdot \x_o)\cdot
(\zH^w)^{-2\rho}\cdot (\zH^w)^{2\rho}\cdot
\cK(aw)\ dn \ da \\
&=& \sum_{w\in \cW} \int_A \cR(f)(a\zH^w  \cdot \xi_o)\cdot
(\zH^w)^{2\rho}\cdot \cK( aw ) \ da \\
&=& (-1)^{k-1}i\left( \int_\R \frac{\cR(f)(a_{t+i\frac{\pi }{2}}\cdot \xi_o)}
{1- e^{-(t+ i\frac{\pi }{2})}}\, dt- \int_\R \frac{\cR(f)(a_{t-i\frac{\pi }{2}}\cdot \xi_o)}
{1- e^{-(t- i\frac{\pi }{2})}}\, dt\right) \, .
\end{eqnarray*}

Consider the strip domain $S=\{ z\in \C\mid |\im z|\le\frac{\pi
}{2}\}$. By our assumption on $f$, the function
$$S\ni z\mapsto F(z)=
\frac{i\cR(f)(a_z\cdot \xi_o) }{ 1- e^{-z}}\in \C$$
defines a meromorphic function $F$ on $\Int S$ which 
extends continuously up to the boundary 
and which has at most a simple pole at $z=0$.
Thus the Residue theorem yields that

$$\cC(f)(\xi_o) = (-1)^k  2\pi i \cdot {\rm Res}(F,0)= (-1)^{k-1} 2\pi \cdot \cR(f)(\xi_o) $$
and this concludes the proof of our theorem.
\end{proof}

\begin{remark} (a) We mention that the geometric pairing $\xi\cdot y$ can be
expressed  using the previously defined power functions
$f_j$ (cf.\ \ref{ef}):
$$\xi\cdot y= f_1(\xi^{-1}y)\, .$$
\par\noindent (b) The assumption that $n$ is even is not a real restriction as 
one can slightly modify $\cC$ so that it works for all parities
(see our discussion in the next section). 
\end{remark}

\section{The holomorphic Radon transform as Cauchy integral II:
Cayley type spaces}
\noindent
The technique used for the hyperboloid in the previous can be used 
for NCC-spaces of Cayley type as well. 
Let us recall that Cayley type spaces are those which are
associated to Euclidean Jordan algebras $V$: i.e. $X$ is a tube domain
associated to $V$ and $H$ is the structure group of the
cone of squares in $V$.
In terms of the set of restricted roots $\Sigma$, this means that $\Sigma$ is
of type $C_n$, say

$$\Sigma=\left\{ \frac{1}{2}(\pm \gamma_i \pm \gamma_j)\mid 1\leq i, j\leq n \right\} \bs \{ 0\}\, . $$

We assume now that $Y=G/H$ is of Cayley type. Define
 $T_j\in \fa$ by $\gamma_i(T_j)=\delta_{ij}$, then

\begin{equation} \label{eq=om} \Omega_H=\bigoplus_{j=1}^n \left ]-\frac{\pi }{2}, \frac{\pi }{2}\right[ T_j\, .
\end{equation}

As a basis of $\Sigma$ we shall choose
$$\Pi=\left\{ \frac{1}{2}(\gamma_1-\gamma_2), \ldots, \frac{1}{2}(\gamma_{n-1}-\gamma_n),\gamma_n\right\}\, .$$
Obviously $\omega_1=-\gamma_1$ is a fundamental spherical lowest weight and accordingly
$f_1(g)=\la \pi_1(g)\eta_j, v_1\ra$ defines a holomorphic
function on $G_\C$ with $f_1(g) =a(g)^{-\gamma_1}$ for $g\in N_\C A_\C K_\C$.
The analogue of Lemma \ref{lem=hyp} now reads
as follows.

\begin{lemma} \label{lem=tra}
Let the notation be as above. Then the following holds:
\begin{enumerate}
\item  $f_1(\exp(i\Omega_H)Gz_H)\subseteq \C \bs \R^\times$.
\item For all $0\leq c < 1$ there exists $C>0$ such that
$$
 (\forall g\in G)(\forall Z\in \Omega_H)\quad |1 -f_1(\exp(ic Z)gz_H)| > C\, .$$
\end{enumerate}
\end{lemma}

\begin{proof} First it is clear from (\ref{eq=om}) that
\begin{equation}\label{e11}
\exp(i\Omega_H)^{\omega_1}=\{ z\in \C \mid \mathrm{Re} z>0\}\, .
\end{equation}
Next recall that $\bigcup_{w\in \cW} NA wH$ is dense in $G$ and that
\begin{equation} \label{e22}
f_1(nawh)=a^{\omega_1} z_H^{w\omega_1}\in \R^+\{ -i, i\}\, .
\end{equation}
Combine (\ref{e11}) and (\ref{e22}) and the assertions follow.
\end{proof}
For $0<c<1$ define a $G$-subdomains of $\Xi_+$
by
\begin{equation}\label{eq-Xc}
\Xi_c=G\exp (ic\Omega_H)\cdot\xi_o\, .
\end{equation}
Note that $-\gamma_j$, $j\not= 1$ is not a fundamental spherical
lowest weight. Therefore, for $1\leq j\leq n$ we define a meromorphic function on $G_\C$ directly by
$$h_j(g)=a(g)^{-\gamma_j} \quad\hbox{for $g\in N_\C A_\C K_\C$}\, .$$
Note that $h_1=f_1$ and, in the same manner
as in Lemma \ref{lem=tra}, one establishes that
\begin{equation} h_j(\exp(i\Omega_H) Gz_H)\subseteq \C \bs \R^\times \cup\{\infty\}
\end{equation}
In particular we see that the prescription
$$\cK: \Xi_+ \times Y\to \C , \ \ (\xi, y)\mapsto
\prod_{j=1}^n \frac{1}{ 1-h_j(\xi^{-1}y)}$$
defines an analytic function,  $CR$ in the
first variable and bounded on all subsets $\Xi_c\times Y$.

\begin{remark} Alternatively, the kernel $\cK$ can be
expressed in Jordan algebra terms.
Let $V$ be the Euclidean Jordan algebra
associated to $X$ and $W\subset V$ its
cone of squares. Form the tube domains
$\cT^\pm = V \pm iW \subset V_\C$.
Then $X=\cT^+$ and $D\simeq \cT^+\times \cT^-$ with
$X$ realized in $D$ via the map
$x\mapsto (x,\oline x)$. Write
$\Delta_j$ for the power functions on $V_\C$ (i.e.
generalized principal minors).
Then, on $D$, one has
$$h_j(z,w)=\frac{\Delta_j(z-w)}{\Delta_{j-1}(z-w)} \qquad\hbox{for $(z,w)\in D$}$$
with the understanding that $\Delta_0 \equiv 1$
Thus $\cK$, when considered as a function on $D$, is given by
$$\cK(z,w)=\frac{\Delta_1 (z-w)\cdot \ldots\cdot \Delta_{n-1}(z-w)}{ \prod_{j=1}^n (\Delta_{j-1}(z-w) -\Delta_j(z-w))}
\, .$$
\end{remark}

However, $\cK$ is not our Cauchy-kernel yet
and some small modification is needed.
For that recall the decomposition
$\Xi_+\simeq G/MN \times \exp(i\Omega_H)$.
For a unitary character $\chi\in \widehat{ T/F}$, where
$F=K\cap T$ is the canonical finite 2-group as usually,  define
the space of {\it $\chi$-twisted $CR$-functions} by
\begin{equation}\label{eq-twisted}
CR_\chi(\Xi_+)=\{ f(ga\cdot \xi_o)=h(ga\cdot\xi_o) \chi^{-1}(a)\mid h\in CR(\Xi_+)\}\, .
\end{equation}
Let $\gamma_0=\gamma_1+\ldots + \gamma_n$. In the sequel we will fix $\chi$ to be
$$\chi= -2\rho +\gamma_0$$
and notice that $\chi={\bf 1}$ for $G=\Sl(2,\R)$.
For $g\in N_\C AT K_\C$ write $t(g)\in T/F$  for the
compact middle part of $g$ and
define the  Cauchy-kernel $\cK$ by
\begin{equation}\label{eq-CKer}
\cK_\chi (\xi, y) =\cK(\xi,y) \chi(t(\xi^{-1}y))\, .
\end{equation}
and notice that $\cK_\chi$ is defined whenever $\xi^{-1}y\in N_\C A_\C K_\C$
and, when defined, is in $CR_\chi$ as a function of the
first variable. We note that $\cH_\chi$ is $G$-invariant and hence
corresponds to a function of one variable $\cK_\chi (g)=\cK_\chi (\xi_o,g\cdot y_0)$.
As before, $\cK_\chi$ is $N\times H$-invariant, and will be identified with
left $N$-invariant function on $Y$.

We
define the \textit{twisted
Cauchy transform} by
$$\cC_\chi: \cS(Y)\to CR_\chi(\Xi_+), \ \ \cC_\chi(f)(\xi)=\int_Y f(y) \cK_\chi(\xi, y) \ dy$$
and the {\it twisted holomorphic Radon transform} by
$$\cR_\chi : \cH^2(D)_{00}\to CR_\chi (\Xi_+), \ \ f\mapsto
\left( \xi=ga\cdot\xi_o\mapsto a^{-\gamma_0}
\int_N f(gna\cdot x_o)\ dn\right)\, .$$
We come to the main result of this section.

\begin{theorem} Suppose that $Y=G/H$ is of Cayley type. Let $f\in \cH^2(D)_{00}$. Then

$$\cC_\chi (f)(\xi)=(-2\pi i )^n \cdot \cR_\chi (f)(\xi) \qquad (\xi\in \Xi_+)\, .$$
\end{theorem}

\begin{proof} Both $\cC_\chi (f)$ and $\cR_\chi (f)$ are $CR$-functions, and so it
is sufficient to show that both coincide on $G/MN\subset \Xi$.
Next,  by $G$-equivariance of both maps, it is enough  to show that
\begin{equation}\label {e111}  \cC_\chi (f)(\xi_o)=(2\pi)^n \cdot \cR_\chi(f)(\xi_o)\, .\end{equation}

We now get, as in the proof of Theorem 4.2, 
 for the left hand side:

\begin{eqnarray*} \cC_\chi(f)(\xi_o) &=& \int_Y  f(y) \cdot \cK_\chi( y)\ dy \\
&=& \ldots \\
&=& \sum_{w\in \cW/ \cW_H} \int_A \cR_\chi(f)(a\zH^w  \cdot \xi_o)\cdot (\zH^w)^{\gamma_0}\cdot \cK_\chi( a\zH^w\cdot x_o) \ da \\
&=& \sum_{w\in \cW/ \cW_H} \int_A \cR(a\zH^w  \cdot \xi_o)\cdot (\zH^w)^{\gamma_0}\cdot \cK(a\zH^w\cdot x_o) \ da \, .
\end{eqnarray*}
Specifically we have $\cW/ \cW_H\simeq (\Z_2)^n$ and $\zH^{w\gamma_0}=i^n {\rm sgn}(w)$

\begin{eqnarray*} \cC_\chi(f)(\xi_o) &=&
i^n \sum_{\eps\in (\Z_2)^n} {\rm sgn}(\eps)
 \int_{\R^n} \cR(f)(\exp\left(\sum_{j=1}^n (t_j+i\eps_j \frac{\pi }{2})T_j\right)\cdot \xi_o)\cdot\\
& & \quad \cdot  \cK( \exp\left(\sum_{j=1}^n (t_j+i\eps_j
 \frac{\pi }{2})T_j\right)\cdot x_o) \ dt
\end{eqnarray*}
Next observe that
$$\cK( \exp\left(\sum_{j=1}^n (t_j+i\eps_j \frac{\pi }{2})T_j\right)\cdot x_o)=\prod_{j=1}^n \frac{1}{ 1- e^{-(t_j+ i \eps_j \frac{\pi }{2})}}\, .$$

Let us consider the multi strip domain $S=\{ z\in \C^n \mid |\im
z_j|\le \frac{\pi }{2}\}$. By our assumption on $f$, the prescription
$$S\ni z\mapsto \cR(f)(\exp\left(\sum_{j=1}^n z_j T_j\right)\cdot  \xi_o) \cdot \prod_{j=1}^n
\frac{1}{ 1- e^{-z_j}}\in\C$$
defines a meromorphic function $F$ on $\Int S$ which 
extends continuously up to the boundary and with 
at most a simple multi-pole at $z_j=0$.
Thus iteratively applying the the Residue theorem yields that

$$\cC_\chi (f)(\xi_o) = (-2\pi i)^n \cdot {\rm Res}(F,0)= (-2\pi i )^n \cdot \cR(f)(\xi_o)
=(-2\pi i)^n \cR_\chi(f)(\xi_o) $$
and it concludes the proof of our theorem.
\end{proof}

\section{Some remarks on the inversion of the holomorphic Radon transform}
\noindent
The inversion of the real horospherical
transform on $X$ can be analytically continued
to give the inversion of the holomorphic
horospherical transform. The dual
transform is given by integration over the real form
$S_\R^K$ of $S(z)$. However, there is a second non-compact
real form $S_\R^H(z)$ of $S(z)$ which gives rise to a different
dual transform and inversion. This is the topic of this section.
Mainly we will focus on $G=\Sl(2,\R)$.

\par We begin with the definition of an appropriate function
space. Let us denote by  $\cF(\Xi_+)$ denote the space
of $CR$-functions on $\Xi_+$ which extend continuously
to $G\exp(i\overline \Omega_H)\cdot\xi_o$
such that $H\ni h \mapsto f(ghz_H\cdot \xi_o)\in\C$ is
integrable for all $g\in G$.
For those functions we define the dual
Radon transform
by
$$\cF(\Xi_+)\to C(Y), \ \ \phi\mapsto \phi^\vee$$
with

$$\phi^\vee(y) =\int_H f(ghz_H\cdot\xi_o)\ dh
=\int_{S_\R(y)} f   \qquad (y=g\cdot y_o\in Y)\, .$$
Clearly, this is a $G$-equivariant mapping.

We would like to understand the relation between $\cR$ and $\phi\mapsto \phi^\vee$.
In this context we would like to mention the
result  in \cite{gko2} for the  holomorphic discrete series:
there exists a differential operator ${\mathcal L}$ such that
$({\mathcal L}\cR (f))^\vee =f $. Hence it is natural to
ask whether a similar statement would hold true for the
most continuous spectrum considered in this paper.
It will turn out that the situation different for the
most continuous series in the sense that
the inverting operator ${\mathcal L}$ is not a
differential operator. We give a detailed  discussion for the
the basic example.

\subsection{The example of $G=\Sl(2,\R)$}
For this paragraph we let $G=\Sl(2,\R)$ with the usual choices
$$A=\left\{\begin{pmatrix} t & 0 \\ 0 & \frac{1}{ t}\end{pmatrix}\mid t>0\right \},
\quad N=\left\{\begin{pmatrix} 1 & x \\ 0 & 1\end{pmatrix}\mid x\in \R\right \},
$$
and $K=\SO(2,\R)$. Let $f\in \cH^2(D)_{00}$ be a $K$-invariant
function. In the sequel we will identify $\fa_\C^*$ with $\C$
via the assignment
$$\C\ni  \lambda\mapsto \lambda\cdot\rho\in \fa_\C^*\,. $$
In this coordinates one has
$$\bc(\lambda)= \pi^{-{1/ 2}} \frac{\Gamma (\lambda/2)}{ \Gamma((\lambda+1)/2)}
\quad \hbox{and} \qquad |\bc(\lambda)|^{-2}= \frac{i\pi\lambda}{ 2} \tanh
\left(\frac{i\pi\lambda}{2}\right)$$
We know that $f|_X\in L^2(X)$ and,  as $f$ is $K$-invariant, we can write
$$f(x)=\frac{1}{2} \int_\R \hat f(i\lambda)\phi_{i\lambda}(x)\,
\frac{d\lambda}{|\bc(i\lambda)|^2} \qquad (x\in X)\, .$$
Applying $\cR$ yields that

$$\cR(f)(\xi)=\frac{1}{2} \int_\R \hat f(i\lambda) a(\xi^{-1})^{\rho(1+i\lambda)} \
 d\lambda\, , $$
and thus

$$\cR(f)^\vee (y_o)=\frac{1}{2} \int_H
\int_\R \hat f(i\lambda) a(\zH^{-1}h)^{\rho(1+i\lambda)} \ d\lambda \ dh \, .$$
Now, for $h=\begin{pmatrix} \cosh t & \sinh t \\
 \sinh t & \cosh t\end{pmatrix}\in H=\SO_e(1,1)$ one has
$\zH^{-1} h =  \begin{pmatrix} e^{-i\frac{\pi}{ 4}}
\cosh t & e^{-i\frac{\pi}{ 4}}\sinh t \\ e^{i\frac{\pi }{4}}\sinh t &
e^{i\frac{\pi }{4}}\cosh t\end{pmatrix}$ and so
$$a(\zH^{-1} h)^\rho=\left(\frac{1}{i(\sinh ^2 t + \cosh ^2 t)}\right)^\frac{1}{2}
=e^{-i\frac{\pi }{4}} \cdot \frac{1}{ (\cosh 2t)^\frac{1}{2}}$$
Therefore we obtain that

\begin{equation}\label{I1} \cR(f)^\vee (y_o)=\frac{1}{2} \int_\R  \int_\R \hat f(i\lambda)
\cdot \frac{e^{\frac{\pi}{4}(\lambda -i)} }{ (\cosh 2t)^{\frac{1}{2}(1+i\lambda)}}
 \ d\lambda\ dt \, .\end{equation}

\begin{lemma}
$$\int_\R  \frac{1}{ (\cosh 2t)^{\frac{1}{2}(1+i\lambda)}}\ dt
 =\frac{1}{2}B(1/2, (1+i\lambda)/4)=\frac{1}{2} \frac{\Gamma((1+i\lambda)/4) \Gamma(1/2)}{ \Gamma((3 + i\lambda)/4)}$$
\end{lemma}

\begin{proof}  Let us denote the integral on the left hand side by $I(\lambda)$.
With the substitution $u=\cosh 2 t$ we obtain
\begin{eqnarray*}
I(\lambda) & =& \int_1^\infty \frac{1}{ u^{\frac{1}{2}(1+i\lambda)}}
\frac{1}{ (u^2 -1)^{1/2}} \ du \\
&=&\frac{1}{2} \int_1^\infty v^{-\frac{1+i\lambda}{ 4}}
 (v -1)^{-{1/2}} v ^{-1/2}\ dv\qquad(v=u^2)\\
&=&\frac{1}{2}B(1/2, (1+i\lambda)/4)
\end{eqnarray*}
as
$B(p,q)=\int_1^\infty u^{-(p+q)}(u-1)^{p-1}\, du$.
\end{proof}
Using this, we get:
\begin{equation}\label{I2} \cR(f)^\vee (y_o)=\frac{1}{4} \int_\R \hat f(i\lambda)
\cdot e^{\frac{\pi}{4}(\lambda -i)} B(1/2, (1+i\lambda)/4)
 \ d\lambda\, .\end{equation}
Define
$$C_1(\lambda)= e^{\frac{\pi}{4}(\lambda -i)} B(1/2, (1+i\lambda)/4)
+  e^{-\frac{\pi}{4}(\lambda +i)}B(1/2, (1-i\lambda)/4) $$
and note that $\lambda\mapsto \hat f(i\lambda)$ is an even function. Thus
(\ref{I2}) yields that

\begin{equation}\label{I3} \cR(f)^\vee (y_o)=\frac{1}{2} \int_\R
\hat f(i\lambda) \cdot C_1(\lambda)
 \ d\lambda\, .\end{equation}
 By the Fourier inversion formula, we have
\begin{equation} \label{I4} f(y_o)=\frac{1}{2} \int_\R \hat f(i\lambda)
\phi_{i\lambda}(y_o) {d\lambda\over |\bc(\lambda)|^2}\, .
\end{equation}
Now, from   the special values of the Gau\ss{} hypergeometric function
we get

$$\phi_{i\lambda}(y_o)=F(1/4+i\lambda/4, 1/4 - i\lambda/4, 1; 1)=
\frac{1}{2B((3-i\lambda )/4,(3+i\lambda )/4)}\, .$$
We therefore define
$$C_2(\lambda) =\frac{1/2}{B((3-i\lambda )/4,(3+i\lambda )/4)\, |\bc(i\lambda)|^2}$$
and note that (\ref{I4}) transforms into
\begin{equation} \label{I5} f(y_o)=\frac{1}{2} \int_\R \hat f(i\lambda)
\cdot C_2(\lambda) \ d\lambda
\end{equation}
In the next step we want to compare the expressions $C_1(\lambda)$ and $C_2(\lambda)$.
If there would exist a differential operator ${\mathcal L}$, then there should be
a polynomial $g(\lambda) $ such that $g \cdot C_1=C_2$.

We first consider

$$c_1(\lambda):=C_1(\lambda) \cdot \left({\Gamma(1/2)\over \Gamma(3/4 - i\lambda/4)
\Gamma(3/4 + i \lambda/4)}\right)^{-1}=C_1^+(\lambda) +
C_1^-(\lambda)$$
with
$$C_1^\pm (\lambda) =e^{\pm {\pi\lambda \over 4}-i\frac{\pi }{4}}\cdot
{\Gamma((1\pm i\lambda)/4) \Gamma(1/2)\over \Gamma(3/4 \pm i\lambda/4)}\cdot
\left( {\Gamma(1/2)\over \Gamma(3/4 - i\lambda/4)\Gamma(3/4 + i \lambda/4)}\right)^{-1}\, .
$$
We focus on $C_1^+$ and obtain

\begin{eqnarray*} C_1^+ (\lambda) & = &
e^{\frac{\pi}{4}(\lambda -i)}\cdot
{\Gamma((1+ i\lambda)/4) \Gamma(1/2)\over \Gamma(3/4 + i\lambda/4)} \cdot\left(
{\Gamma(1/2)\over \Gamma(3/4 - i\lambda/4)\Gamma(3/4 + i \lambda/4)}\right)^{-1}
\\
& =& e^{\frac{\pi}{4}(\lambda -i)}\cdot
\Gamma((1+ i\lambda)/4) \Gamma(3/4 - i\lambda/4)\\
& =& e^{\frac{\pi}{4}(\lambda -i)}\cdot
\Gamma((1+ i\lambda)/4) \Gamma(1 + (-1/4 - i\lambda/4))\\
& =& e^{\frac{\pi}{4}(\lambda -i)}\cdot
(-1/4 - i\lambda/4)\Gamma((1+ i\lambda)/4) \Gamma(-(1 +i\lambda)/4)\\
& =& {e^{\frac{\pi}{4}(\lambda -i)}\pi \over \sin \pi (1/4 +i\lambda/4)}
\end{eqnarray*}
Likewise we obtain

$$C_1^-(\lambda)=C_1^+(-\lambda)
 = {e^{-\frac{\pi}{4}(\lambda +i)}\pi \over \sin \pi (1/4 -i\lambda/4)}$$
and so

\begin{eqnarray*} c_1(\lambda) &=& \frac{e^{\frac{\pi}{4}(\lambda -i)}\pi
}{\sin \pi (1/4 +i\lambda/4)} +
 \frac{e^{-\frac{\pi}{4}(\lambda +i)}\pi }{ \sin \pi (1/4 -i\lambda/4)}\\
&=& \frac{\pi}{ i} \frac{e^{\frac{\pi}{4}(\lambda -i)}}{ \sinh
 \pi (\lambda/4- i/4)}
+  \frac{\pi}{i} \frac{e^{-\frac{\pi}{4}(\lambda +i)} }{ \sinh \pi
(-\lambda/4- i/4)}\\
&=& \frac{\pi}{i} \frac{\cosh ({\pi\lambda }{ 4}-i\frac{\pi }{4}) +
\sinh({\pi\lambda \over 4}-i\frac{\pi }{4})}{ \sinh \pi (\lambda/4- i/4)}
+  {\pi\over i}{\cosh (-{\pi\lambda \over 4}-i\frac{\pi }{4}) +
\sinh(-{\pi\lambda \over 4}-i\frac{\pi }{4})
\over \sinh \pi (-\lambda/4- i/4)}\\
&=& {\pi\over i} \cdot\left( 2 + \coth (\frac{\pi}{4}(\lambda -
i)) + \coth (-\frac{\pi}{4}(\lambda+i))\right)
\end{eqnarray*}
Now define $g(\lambda)$ by the requirement
$$g(\lambda) c_1(\lambda)= {1\over |\bc(i\lambda)|^2}=
\frac{\pi\lambda}{ 2} \tanh\left(\frac{\pi\lambda}{ 2}\right)$$
Now note that

\begin{eqnarray*}c_1(\lambda + i) & =&
 {\pi\over i} \cdot\left( 2 + \coth ({\pi\lambda \over 4}) + \coth (-{\pi\lambda \over 4}-i\frac{\pi }{2})\right)\\
 &=& {\pi\over i} \cdot\left( 2 + \coth ({\pi\lambda \over 4}) - \tanh ({\pi\lambda \over 4})\right)\\
&=& {\pi\over i} \cdot\left( 2 + {2\over \sinh ({\pi\lambda \over 2})}\right)
\end{eqnarray*}
and thus we get
$$g(\lambda+i ) {\pi\over i} \cdot\left( 2 + {2\over \sinh ({\pi\lambda \over 2})}\right)
= {\pi(\lambda+i)\over 2} \coth\left({\pi\lambda\over 2}\right)\, .$$
Further manipulation then yields that
$$g(\lambda)= {i\lambda\over 4} \cdot {\sinh({\pi\lambda\over 2})\over 1- \cosh ({\pi\lambda\over 2})}$$
and it is obvious that $g$ is not a polynomial function.
Since
$$g(\lambda) C_1(\lambda)= C_2(\lambda) $$
it is now clear that there exists {\it  no} differential operator ${\mathcal L}$ which inverts the holomorphic Radon transform.
However the function $g(\lambda)$ defines us a spectral multiplier which is
a  pseudo-differential operator which we now call ${\mathcal L}$.
We summarize our discussion.

\begin{theorem} Let $G=\Sl(2,\R)$ and ${\mathcal L}$ the spectral
multiplier defined by the function $g(\lambda)= {i\lambda\over 4} \cdot {\sinh({\pi\lambda\over 2})\over 1- \cosh ({\pi\lambda\over 2})}$.
Then for $f\in \cH^2(D)_0$ such that $\cR(f)\in \cF(\Xi_+)$ one has
$$f = ({\mathcal L} \cR (f))^\vee\, .$$
\end{theorem}

\section{Geometric definition of the Hardy space}
\noindent
This final section deals with the structure
of the
Hardy space $\cH^2(D)$. It allows independent
reading and is of independent interest.
\par Initially, the Hardy space was defined
spectrally (see \cite{GKO}). Below we will show how
to define the Hardy space geometrically, i.e. we give
geometric definition of the norm on $\|\cdot\|_H$ on
$\cH^2(D)$ through
$G$-orbit integrals on $D$. For that we start
by recalling the orbital integral $\bO_h$ and
the pseudo-differential operator $\cD$ introduced
also used in \cite{kos05}.

\par In this section $Y=G/H$ can  be an arbitrary NCC-space. 

\subsection{$G$-orbit integrals on the domain $D$}
For a sufficiently decaying functions $h$ on
$D$ we define its {\it $G$-orbit integral} on $D$
as the following function on $i2\Omega_H$

$$\bO_{h}(iX)=\int_G h(g\exp(i\frac{1}{ 2}X)\cdot x_o) \ dg\qquad (X\in 2\Omega_H)\ .$$
For $f\in \cH^2(D)$ we notice that $|f|^2$ is a sufficiently decaying
function on $D$, i.e. $\bO_{|f|^2}(iX)$ is finite for all
$X\in 2\Omega_H$. Moreover,  in view of (\ref{eq=Gutzmer})
we see that $\bO_{|f|^2}$ has a natural holomorphic extension
to a holomorphic function on the abelian tube domain
$\cT(2\Omega_H)=\fa+i2\Omega_H$, namely
\begin{equation}\label{eq=exgutz}
\bO_{|f|^2}(Z)=\int_{\mathcal X} |\hat f(b,\lambda)|^2
\varphi_\lambda(\exp(Z)) \ d\mu_{\mathcal X}(b,\lambda)
\qquad (Z\in \cT(2\Omega_H))
\end{equation}

\subsection{A certain  pseudo-differential operator}
Define a space $\cF(\cT(2\Omega_H))$ of $\cW$-invariant holomorphic functions on
the tube domain $\cT(2\Omega_H)$ by the following property:
$f\in\cF(\cT(2\Omega_H))$ if $f$ can be written  as
$$f(Z)=\int_{i\fa_+^*} h(\lambda)\varphi_\lambda(\exp(Z))
\frac{d\lambda}{ |\bc(\lambda)|^2}\qquad (Z\in \cT(2\Omega_H))$$
where $h\in L^1(i\fa_+^*, \frac{\COSH(\lambda)}{ |\bc(\lambda)|^2}d\lambda)$.
If $Q\subset \cT(2\Omega_H)$ is compact, then there exists a
constant $C_Q>0$ such that
$$(\forall \lambda\in i\fa^*)\qquad \sup_{X\in Q}
|\varphi_\lambda(\exp(i2X))| \le C_Q \COSH(\lambda)\, .$$
As $\frac{1}{|\bc(\lambda)|^2}$ is at most of polynomial
growth, it follows that $f$ is indeed holomorphic  and $\cW$-invariant. Moreover,
$f$ is uniquely determined by $h$. It follows
from our discussion that the prescription

$$\cD: \cF(\cT(2\Omega_H))\to \cO (\cT(2\Omega_H))^\cW;\   (\cD F)(Z)=\int_{i\fa_+^*}
h(\lambda) \sum_{w\in \cW} e^{\lambda(wZ)}
\frac{d\lambda}{ |\bc(\lambda)|^2}$$
is a well defined linear mapping.

\begin{remark}The operator $\cD$ is a pseudo-differential  operator and
a differential operator if all root multiplicities are even. The operator $\cD$
is related to the Abel transform as explained in
\cite{kos05}, Remark 3.2.
\end{remark}

\begin{example} In this example we discuss the operator $\cD$
when the underlying group $G$ is complex. Then 
$\cD$ is a differential operator of a particularly
nice form.

\par   If $G$ is complex, then there is an explicit formula
for spherical functions, due to Harish-Chandra:

$$\varphi_\lambda(\exp(Z))={\bc}(\lambda)
\frac{\sum_{w\in \cW} \eps(w)
e^{\lambda(wZ)}}{\prod_{\alpha\in \Sigma^+} 2\sinh\alpha(Z)}$$
for all $Z\in \cT(2\Omega_H)$. The $\bc$-function has the familiar form

$$\bc(\lambda)=\frac{\prod_{\alpha\in \Sigma^+} \la \rho, \alpha\ra}{
\prod_{\alpha\in \Sigma^+} \la \lambda,\alpha\ra}\ .$$
For each $\alpha\in \Sigma$ let $A_\alpha\in \fa$ be such that
$\alpha=\la\cdot, A_\alpha\ra$. Furthermore let $ \partial_\alpha$
be the partial derivative on $\cT(2\Omega_H)$ in direction $A_\alpha$.
Define a partial differential operator on $\cT(2\Omega_H)$ by
$ \partial_{\Sigma^+}=\prod_{\alpha\in \Sigma^+} \partial_\alpha$.
Finally with $J(Z)=\prod_{\alpha\in \Sigma^+} 2\sinh\alpha(Z)$ we
declare a differential operator on $\cT(2\Omega_H)$ by
$$\cD=\mathrm{const}\cdot
\partial_{\Sigma^+}\circ J\ .$$
with $\mathrm{const}= \prod_{\alpha\in \Sigma^+} \la \rho, \alpha\ra $.
The relation
$$\cD(\varphi_\lambda\circ \exp)(Z)=\sum_{w\in \cW}  e^{\lambda(wZ)}$$
is now obvious.
\end{example}

\subsection{The geometric norm}
For a function $f\in \cH^2(D)$ let us write
$\|f\|_H$ for its norm as before.
By Lemma \ref{eq=hn} this norm
is given by

\begin{equation}\label{eq=spectral}
\|f\|_H = \int_{{\mathcal X}} |\hat f(b,\lambda)|^2
\COSH(\lambda) \ d\mu_{\mathcal X}(b,\lambda)
\end{equation}
The objective of this section is to express $\|f\|_H$ in terms
of the much more geometric orbital integrals $\bO_{|f|^2}$.
Our result is as follows.

\begin{theorem}\label{th=Hardy} Let $f\in \cH^2(D)$. Then
$\bO_{|f|^2}\in \cF(\cT(2\Omega_H))$ and the Hardy space
norm $\|f\|_H$ of $f$ is given by

\begin{equation} \label{eq=norm}
\|f\|_H =\sup_{X\in 2\Omega_H} \frac{\cD\bO_{|f|^2}(iX)}{ |\cW_H|} .
\end{equation}
In particular, the Hardy space $\cH^2(D)$ can be defined as
\begin{equation}\label{eq=char}
\cH^2(D)=\{ f\in \cO(D)\mid \bO_{|f|^2}\in \cF(\cT(2\Omega_H))\
\sup_{X\in 2\Omega_H} \frac{\cD\bO_{|f|^2}(iX)}{ |\cW_H|}<\infty\}\ .
\end{equation}
\end{theorem}

\begin{proof} Fix $f\in \cH^2(D)$. By equation (\ref{eq=exgutz})
we have
$$\bO_{|f|^2}(Z)=\int_{{\mathcal X}} |\hat f(b,\lambda)|^2
\varphi_\lambda(\exp(Z))\ d\mu_{{\mathcal X}}(b,\lambda)$$
for all $Z\in \cT(2\Omega_H)$.
By  the spectral definition of $\cH^2(D)$ it follows that
$$h(\lambda):=\int_B |\hat f(b,\lambda)|^2 \ db
$$
defines  a function $h\in L^1(i\fa_+^*,
\frac{\COSH(\lambda)}{ |\bc(\lambda)|^2}d\lambda)$.
Thus $\bO_{|f|^2}\in \cF(\cT(2\Omega_H))$ and the application of $\cD$  to  $\bO_{|f|^2}$ yields
$$ (\cD\bO_{|f|^2})(Z)=\int_{\mathcal X} |\hat f(b,\lambda)|^2
 \sum_{w\in \cW}e^{\lambda(Z)}\ d\mu_{{\mathcal X}}(b,\lambda)\, .$$
 Now notice that
$$\frac{1}{ |\cW_H|} \sum_{w\in \cW}e^{\lambda(iwX)}\leq \COSH(\lambda)\qquad (X\in
2\Omega_H)$$
and
$$\sup_{X\in 2\Omega_H} \sum_{w\in \cW}e^{\lambda(iwX)}
=\lim_{X\to \ZH}\sum_{w\in \cW}e^{\lambda(i2wX)}
=|\cW_H|\cdot \COSH(\lambda)\ .$$
The claim follows now from the spectral definition of the norm
in $\cH^2(D)$. Finally, backtracking
the steps of the proof readily yields (\ref{eq=char}).
\end{proof}

\begin{remark} Some comments  on the geometric Hardy norm
$$\|f\|_H =\sup_{Z\in \Omega_H} \frac{1}{ |\cW_H|}\cdot(\cD\bO_{|f|^2})(Z)$$
seem to be appropriate. Usually, in the theory of Hardy
spaces (e.g. Hardy space on the upper half plane) one takes
the supremum over a family of $L^2$-integrals over totally real
submanifolds. In our case one takes a supremum over $G$-orbits,
which for the exception of the orbit through the origin, are
never totally real.
  Secondly, we find the appearance of the pseudo differential  operator $\cD$
interesting. In the context of Hardy spaces it might be novel.
\end{remark}

\subsection{The $K$-invariant case}
In this subsection we give another description of the subspace
$\cH^2(D)^K$ using the Abel transform and the results in Appendix A.
We start by noting the following
simple connection between the Abel transform and the Fourier transform of
a $K$-invariant function. For that we note first, that $\hat{f}(b,\lambda)$ is independent
of $b\in B$ if $f$ is $K$-invariant. We write then simply $\hat{f}(\lambda)$ and
note that $\hat{f}$ is $\cW$-invariant. Furthermore
\begin{eqnarray*}
\hat{f}(\lambda )&=&\int_X f(x)a(x)^{\rho - \lambda }\, dx\\
&=&\int_A\int_N f(na  \cdot x_o)a^{-\rho -\lambda}\, dn da\\
&=&\cF_A(\cA (f))(\lambda )
\end{eqnarray*}
where $\cF_A$ stands for the Fourier transform on the abelian group $A$.
Recall, that $\cF : L^2(A)\to L^2(i\fa^*,|\cW |^{-1} d\lambda )$ is a unitary
isomorphism.
Define a multiplication operator $D_\fa$ on $i\fa^*$ by
$D_\fa (F) = \bc (-\lambda  )^{-1}F$ and
denote the corresponding multiplier operator on $A$ by
$D_A$. Let $\Lambda := D_A\circ \cA$. Finally, we define a
multiplier $m$ on $\cW\times i\fa^*$ by $m(s,\lambda )=\bc (-s\lambda )/\bc (-\lambda )$.
We denote by $\tau$ the corresponding representation
$\tau (s)f(\lambda )=m(s^{-1},\lambda )f(s^{-1}\lambda )$.
Then, cf. \cite{os05}, Section 1, in particular Lemma 1.4, we
have a commutative diagram, where each of the maps is
an unitary isomorphism:

\begin{equation}\label{eq=isometries}
 \xymatrix { L^2(A,|\cW|^{-1}da)^{\cW} \ar[r]^\Lambda\ar[d]^{\cF_X} &
 L^2(A,|\cW|^{-1}da)^{\tau (\cW )}\ar[d]^{\cF_A}\\
 L^2(i\frak{a}^*,\frac{d\lambda }{|\cW||\bc (\lambda )|^2})\ar[r]_{D_{\fa}} &
 L^2(i\fa^* ,|\cW|^{-1} d\lambda )^{\tau (\cW)}\,. }
\end{equation}

Recall the Hardy space $\cH^2(T(\Omega_H))$ from Appendix A and its
spectral description in Theorem \ref{co=A1}. It follows then from
Theorem \ref{eq=hn}, and the obvious renormalization
of measures, as we have not included the $2\pi$ in
the exponential function, that $\Lambda (\cH^2(D)^K)\subseteq \cH^2(T(\Omega_H))^{\tau (\cW )}$.
As $\Lambda :L^2(A,|\cW|^{-1}da)^{\cW} \rightarrow
 L^2(A,|\cW|^{-1}da)^{\tau (\cW )}$ is a unitary isomorphism, we get:
\begin{theorem}\label{th-isohardy}
The map $\Lambda :\cH^2(D)^K\to \cH^2(T(\Omega_H))^{\tau (\cW)}$ is a unitary
isomorphism.
\end{theorem}

\begin{example} If $G$ has complex structure, then 
the map $\Lambda$ is a multiplication operator given by
$$\Lambda (f)(a)=\left(\prod_{\alpha\in\Sigma^+}\sinh \la \alpha , \log (a)\ra \right)f(a)\, .$$
\end{example}
We now determine the reproducing kernel for $\cH^2 (D)^K$. One could easily deduct
that from Theorem \ref{th-repkern} but we will give another proof, that
follows similar arguments.
\begin{theorem} The reproducing kernel $\cK (z,w)$ for  $\cH^2 (D)^K$ is given
by
$$\cK (z,w)=\int_{i\fa_+^*}\varphi_\lambda (z)\varphi_{-\lambda }(\overline{w})\, d\mu (\lambda )\, .$$
\end{theorem}

\begin{proof} Let $f\in \cH^2 (D)^K$ and $w\in D$. Recall, that  by
Lemma \ref{eq=hn} we have $\Phi (g)(\lambda )=\hat{g}(\lambda )\COSH (\lambda )$
for all $g\in \cH^2 (D)^K$. Therefore
\begin{eqnarray*}
f(w)&=&\la f,K_w\ra\\
&=&\int_{i\fa_+}\hat{f}(\lambda )\COSH (\lambda )\overline{\widehat{K_w}(\lambda )\COSH (\lambda )}
\, d\mu (\lambda )\\
&=& \int_{i\fa_+^*}\hat{f}(\lambda )\overline{\widehat{K_w}(\lambda )\COSH (\lambda )}
\, \frac{d\lambda}{|\bc (\lambda )|^2}\\
&=&\int_{i\fa_+^*}\hat{f}(\lambda )\varphi_\lambda (w)\, \frac{d\lambda}{|\bc (\lambda )|^2}\, .
\end{eqnarray*}
Thus
$$\widehat{K_w}(\lambda )=
\frac{\varphi_{-\lambda}(\overline{w})}{\COSH (\lambda )}\, .$$
{}From this we now get:
\begin{eqnarray*}
K(z,w)&=&\la K_w , K_z\ra\\
&=& \int_{i\fa_+^*} \widehat{K_w}(\lambda )\COSH (\lambda )
\overline{\widehat{K_z}(\lambda )\COSH (\lambda )}\, d\mu (\lambda)\\
&=&\int_{i\fa_+^*} \varphi_{-\lambda}(\overline{w})
\varphi_{\lambda}(z) \, d\mu (\lambda )
\end{eqnarray*}
and the claim follows.
\end{proof}

\section*{Appendix: Hardy spaces on strip domain}
\noindent
We let $V$ be an Euclidean vector space, e.g. $V=\R^n$ endowed with the
standard inner product. Denote by ${\rm O}(V)$ the orthogonal group of $V$
and let  $\cW\subset {\rm O}(V)$ be a finite subgroup which acts irreducibly
on $V$. We fix $y_o\in V$, $y_o\neq 0$ and set
$$\Omega={\rm int}\left( \{ \hbox{convex hull of }\ \cW(y_o)\}\right)\ .$$
Notice that $\Omega$ is the interior of a compact polyhedron and that
$0\in \Omega$. Write $V_\C :=V\otimes_\R \C\simeq V+iV$
for the complexification of $V$ and define
a tube domain in $V_\C$ by

$$T(\Omega)=V+i\Omega\ .$$
Let us denote by $dx$ the measure $(2\pi )^{-\dim V /2}$ times the
normalized Lebesgue measure on $V$. Then the Fourier transform
$$f\mapsto \cF f(\lambda )= \int_V f(x)e^{-\la \lambda , x\ra}\, dx$$
is
a unitary $L^2$-isomorphism.
$V^*$ is the dual of $V$, and  $\la \lambda  \cdot x  \ra=\lambda (x)$ denotes
the standard duality between $V$ and $V^*$.
 Denote by $\cO (T(\gO ))$ the
space of holomorphic functions on $T(\Omega )$.
The Hardy space $\cH^2(T (\Omega )$ is defined by:
$$\cH^2(T(\Omega)):=\{f\in \cO(T(\Omega))\mid
\|f\|^2_\cH\:=\sup_{y\in \Omega} \int_V |f(x+iy)|^2 \ dx<\infty\}\ .$$

As the Hardy-norm locally dominates the Bergman-norm
on $T(\Omega)$, it follows hence $\cH^2(T(\Omega))$ is complete, i.e.
a Banach space. In fact, $\cH^2(T(\Omega))$ is a Hilbert space
as we will show in a moment.
Then for $f\in \cH^2(T(\Omega))$ and $y\in \Omega$ one has

$$\int_V |f(x+iy)|^2 \ dx =\int_{V^*} |\cF(f|_V)(\xi)|^2 e^{-2 \langle y, \xi\rangle}
\ d\xi$$
which is immediate from \cite{SW} Ch. III, \S2.
It follows that

\begin{equation} \label{eq=A1}
\|f\|^2_\cH =\sup_{y\in \Omega}\int_{V^*} |\cF(f|_V)(\xi)|^2 e^{- 2 \langle y, \xi\rangle}
\ d\xi\ .
\end{equation}
For $y \in V $ define $\COSH,\COSH_y : V^*\to \C$ by
$$\COSH_y (\lambda )=\frac{1}{|W|}\sum_{s\in \cW}
e^{-2 \langle y,s \lambda \rangle}=\frac{1}{|W|}\sum_{s\in \cW}
e^{- 2 \langle sy,\lambda \rangle}$$
and
$$\COSH (\lambda )=\COSH_{y_o}(\lambda )\, .$$

As with $\|\cdot\|_{L^2(V)}$ is $\cW$-invariant and $\cF$ is
$\cW$-equivariant, it follows from
(\ref{eq=A1}) that

\begin{equation} \label{eq=A2}
\|f\|^2=\sup_{y\in \Omega}\int_{V^*} |\cF(f|_V)(\xi)|^2 \COSH(y,\xi)
\ d\xi\ .
\end{equation}

Now, for every $\lambda \in V$, the function
$y\mapsto  \COSH_\y ( \lambda )$ is strictly convex on $\oline \Omega$;
hence we have the inequality
$$(\forall y\in \Omega)(\forall \lambda \in V^*)
\qquad  \COSH_y(\lambda )\leq \COSH (\lambda )\ .$$

and so it follows that

\begin{equation} \label{eq=A3}
\|f\|^2_\cH =\int_{V^*} |\cF(f|_V)(\lambda )|^2 \COSH(\lambda )
\ d\lambda \ .
\end{equation}

Define a $\cW$-invariant measure $\mu_0$ on $V^*$ by
$$d\mu_0 (\lambda )=\COSH (\lambda )d\lambda\, .$$

\begin{theorem}\label{th=A1} The mapping
$$\cH^2(T(\Omega))\to L^2(V^*, d\mu_0),
\ \ f\mapsto \cF(f|_V)$$
is an isometric isomorphism. In particular,
$\cH^2(T(\Omega))$ is a Hilbert space.
\end{theorem}

A (continuous) multiplier on $V^*$ is a continuous map
$m :\cW\times V^*\to \C$ such that for all $s,w\in \cW$ and
$\lambda\in V^*$ we have
$$m(sw,\lambda )=m(s,w\lambda )m(w,\lambda )\, .$$
Assume from now on that $|m(s,\lambda )|=1$ for
all $s\in\cW$ and $\lambda \in V^*$. Then, because of the
$\cW$-invariance of $d\mu_0$, we can define a
unitary representation of $\cW$ on $L^2(V^*, d\mu_0)$ by
\begin{equation}\label{eq-tau}
\tau (s) f(\lambda )=m(s^{-1},\lambda )f(s^{-1}\lambda )\, .
\end{equation}
As the Fourier transform is a unitary isomorphism, we  have
a unitary representation, also denoted by $\tau$,
of $\cW$ on $\cH^2 (T(\Omega ))$ such
that the Fourier transform is an intertwining operator. Denote the space
of $\tau (\cW)$-invariant elements by the superscript $\tau (\cW)$. Then
\begin{corollary}\label{co=A1}
The Fourier transform is a unitary isomorphism
$$\cF :\cH^2(T (\Omega ))^{\tau (\cW)}\to L^2(V^*,d\mu_0 )^{\tau (\cW)}\, .$$
\end{corollary}

The Hardy space $\cH^2(T(\Omega))^{\tau (\cW)}$ being a Hilbert space
of holomorphic functions admits as such a reproducing kernel
function $K(z,w)$ often called the Cauchy-Szeg\"o-kernel.
Notice that for fixed $w\in T(\Omega)$, the function
$K_w(z)$ belongs to $\cH^2(T(\Omega))^{\tau (\cW)}$ and that
$$\langle f, K_w\rangle =f(w) \qquad \hbox{for all $f\in \cH^2(T(\Omega))^{\tau (\cW)}$}\ .$$
Here of course $\la\, \cdot\, ,\, \cdot\, \ra$ denotes the inner product
on $\cH^2(T(\Omega))^{\tau (\cW)}$.

We now determine $K(z,w)$ explicitly.
For that define for $w\in V$, $\COS_w^m:V^*\to \C$ by
$$\COS_{w}^m (\lambda ):=\frac{1}{ |\cW|}
\sum_{s\in \cW} m(s,\lambda )^{-1} e^{ i \la w,s\lambda \ra}$$
and note that
$$\frac{\COS_w^m}{\COSH}\in \cH^2(T(\Omega))^{\tau (\cW)}$$
Write $(\cdot|\cdot)$ for the inner product on
$L^2(V^*,d\mu_0)^{\tau (\cW)}$.
For $f\in \cH^2(T(\Omega))^{\tau (\cW)}$ let $F=\cF (f|_V)$.

It follows from
Corollary \ref{co=A1} that $\la f, K_w\ra =(F|\cF(K_w|_V))$. On the other hand we
have
\begin{eqnarray*}
f(w)&=&\int_{V^*} F(\lambda )e^{- i \la w,\lambda \ra}\ d\lambda\\
&=&\int_{V^*} F(\lambda )
\frac{e^{- i \la w,\lambda \ra}}{ \COSH(\lambda )}\, d\mu_0(\lambda )\\
&=&\int_{V^*} m (s^{-1},\lambda ) F(s^{-1}\lambda )
\frac{e^{- i \la w,\lambda \ra}}{ \COSH(\lambda )}\, d\mu_0(\lambda )\\
&=&\int_{V^*} F(\lambda )  m (s^{-1},s\lambda )
\frac{e^{- i \la w,s \lambda \ra}}{ \COSH(\lambda )}\, d\mu_0(\lambda )\\
&=&\int_{V^*} F(\lambda )\overline{
\frac{\COS_{w}^m(\lambda )}{ \COSH(\lambda )} }\
d\mu_0 (\lambda )
\end{eqnarray*}
and thus
$$\cF (K_w|_V)(\lambda ) = \frac{\COS_{w}^m(\lambda )}{ \COSH(\lambda )}\, .$$

\begin{theorem}\label{th-repkern} The reproducing kernel for the Hardy space $\cH^2(T(\Omega))^{\tau (\cW)}$
is given by
\begin{eqnarray}
K(z,w) &=&\int_{V^*} \frac{\left(\frac{1}{|\cW|}\sum_{s\in \cW}
\frac{
 e^{ i \la s(z),\lambda \ra}}{m(s,\lambda )} \right)\cdot \left(
\oline{\frac{1}{ |\cW|}\sum_{s\in \cW} \frac{
e^{i \la s(w),\lambda \ra}}{m(s,\lambda )}}\right)
}{\COSH(\lambda )} \ d\lambda\label{eq=A4}\\
&=& \int_{V^*} \frac{\COS^m_z(\lambda)}{\COSH (\lambda )}\COS_{-\oline w}^m(\lambda )\, d\lambda\, .
\nonumber
\end{eqnarray}
\end{theorem}\label{ap=K(z,w)}

\begin{example} The equation (\ref{eq=A4}) can be evaluated
in the relevant special cases. Let us for example consider the case
of $V=\R$, $\Omega=]-1,1[$, $\cW={\rm O}(\R)\simeq\{ \one, -\one\}$
and $m(s,\lambda)=1$. Using the standard measure on $\R$
we get the following from (\ref{eq=A4}) and \cite{Bate}, Sect. 1.9, formula (12):

\begin{eqnarray*}
K(z,w)&=& \frac{1}{\sqrt{2\pi}} \int_{\R} \frac{\cos(z\xi) \cos( \oline w\xi )
}{\cosh( 2 \xi )} \ d\xi\\
&=&\sqrt{\frac{\pi }{ 2}}  \frac{\cosh (\frac{\pi}{ 4} z)\cdot\cosh (
\frac{\pi}{ 4} \oline w)}{
\cosh (\frac{\pi}{2} z) + \cosh (\frac{\pi}{ 2} \oline w)}
\, .
\end{eqnarray*}
\end{example}

\end{document}